\documentclass[10pt]{article}
\usepackage{amsfonts,amsmath,latexsym,verbatim,amscd}

\title{Constant mean curvature hypersurfaces \\ 
condensing along a submanifold}

\author{Fethi Mahmoudi \thanks{Email: mahmoudi@univ-paris12.fr}\\ Universit\'e Paris XII
\and
Rafe Mazzeo \thanks{Email: mazzeo@math.stanford.edu. Supported by the NSF
under Grant DMS-0204730}\\ Stanford University
\and
Frank Pacard \thanks{Email: pacard@univ-paris12.fr}\\ Universit\'e Paris XII}

\date{}

\newtheorem{theorem}{Theorem}[section]
\newtheorem{proposition}{Proposition}[section]

\newtheorem{lemma}{Lemma}[section]

\newcommand{\R}{\mathbb{R}}
\newcommand{\N}{\mathbb{N}}

\newcommand{\e}{\varepsilon}
\newcommand{\del}{\partial}
\newcommand{\ds}{}

\newcommand{\calC}{{\mathcal C}}

\newcommand{\calL}{{\mathcal L}}

\newcommand{\calO}{{\mathcal O}}
\newcommand{\calS}{{\mathcal S}}

\begin{document}

\maketitle

\section{Introduction}

Let $S$ be an oriented embedded (or possibly immersed) hypersurface
in a Riemannian manifold $(M^{m+1},g)$. The
shape operator $A_S$ is the symmetric endomorphism of the tangent
bundle of $S$ associated with the second fundamental form of $S$, $b_S$, by
\[
b_S (X, Y) = g_S ( A_S \, X, Y), \quad \forall X,Y \in TS; \qquad
\mbox{here}\qquad g_S = \left. g \right|_{TS}.
\]
The eigenvalues $\kappa_i$ of the shape
operator $A_S$ are the principal curvatures of the hypersurface
$S$. The mean curvature of $S$ is define to be the average of the
principal curvatures of $S$, i.e.
\[
H (S) : = \frac{1}{m} ( \kappa_{1} +\ldots +\kappa_{m}).
\]

Constant mean curvature hypersurfaces constitute a very important
class of submanifolds in a compact Riemannian manifold
$(M^{m+1},g)$. In this paper we are interested in families of such
submanifolds, with mean curvature varying from one member of the
family to another, which `condense' to a submanifold $K^k \subset
M^{m+1}$ of codimension greater than $1$. Under fairly reasonable
geometric assumptions \cite{Maz-Pac}, the existence of such a
family implies that $K$ is minimal. Two cases have been studied
previously: Ye \cite{Ye-1}, \cite{Ye-2} proved the existence of a
local foliation by constant mean curvature hypersurfaces when $K$
is a point (which is required to be a nondegenerate critical point
of the scalar curvature function); more recently, the second and
third authors \cite{Maz-Pac} proved existence of a partial
foliation when $K$ is a nondegenerate geodesic. In this paper we
extend the result and methods of \cite{Maz-Pac} to handle the
general case, when $K$ is an arbitrary nondegenerate minimal
submanifold. No extra curvature hypotheses are required. In
particular, this proves the existence of constant mean curvature
hypersurfaces with nontrivial topology in any Riemannian manifold.

\medskip

Let us describe our result in more detail. Let $K^k$ be a
closed (possibly immersed) submanifold in $M^{m+1}$, $1\leq k\leq m-1$,
and define the geodesic tube of radius $\rho$ about $K$ by
\[
\bar S_\rho :=\{q\in M^{m+1}: \quad \mbox{dist}_g (q,K)= \rho\}.
\]
This is a smooth (immersed) hypersurface provided $\rho$ is
smaller than the radius of curvature of $K$, and we henceforth
always tacitly assume that this is the case. The mean curvature of
this tube satisfies
\[
H ( \bar S_\rho  ) = \frac{n-1}{m} \, \rho^{-1} + \calO (1), \qquad
\mbox{as}\qquad \rho \searrow 0,
\]
with $n = m+1 -k$ and hence it is plausible that we might be able
to perturb this tube to a constant mean curvature hypersurface
with $H \equiv \frac{n-1}{m} \, \rho^{-1}$. This is not quite true
since the mean curvature of $\bar S_\rho $ is not sufficiently
close to being constant, but when $K$ is minimal there is a better
estimate
\[
H ( \bar S_\rho ) = \frac{n-1}{m} \, \rho^{-1} + \calO (\rho),
\]
cf. \S 4. Even in this case, there are other more subtle
obstructions to carrying out this procedure at certain radii
$\rho$ related to eigenvalues of the linearized mean curvature
operator on $\bar S_\rho$, which in turn are related to a genuine
bifurcation phenomenon, at least when $k=1$, \cite{Maz-Pac}. Thus
we do not obtain existence of the constant mean curvature
perturbation for every small radius.
\begin{theorem}
Suppose that $K^k$ is a nondegenerate closed minimal submanifold $1\leq k
\leq m-1$. Then there exists a sequence of disjoint nonempty
intervals $I_i = (\rho_i^-,\rho_i^+)$, $\rho_i^\pm \to 0$, such that
for all $\rho\in I : = \cup_i I_i$, the geodesic tube $\bar S_{\rho}$
may be perturbed to a constant mean curvature hypersurface $S_{\rho}$ with
$H = \frac{n-1}{m } \, \rho^{-1}$. \label{th:existence}
\end{theorem}

The nondegeneracy condition on $K$ is simply that the linearized
mean curvature operator, also called the Jacobi operator, is
invertible; this restriction is quite mild and holds generically
\cite{Whi}. As noted above, this result was already known when
$k=0,1$, but the case $k > 1$ requires a more complicated
analysis. This new approach is inspired by some recent work of
Malchiodi and Montenegro in a somewhat different context
\cite{Mal-Mon-2}, \cite{Mal}.

\medskip

The hypersurface $S_\rho$ is a small perturbation of $\bar S_\rho
$ in the sense that it is the normal graph of some function (with
$L^\infty$ norm bounded by a constant times $\rho^3$) over a
submanifold obtained by `translating' $K$ by a section of its
normal bundle (with $L^\infty$ norm bounded by a constant times
$\rho^2$); we refer to \S 3.1 for the precise formulation of the
construction of $S_\rho$. When $K$ is embedded, then so are the
hypersurfaces $S_{\rho}$ for $\rho$ sufficiently small. In
addition, the hypersurfaces in each of the families
$\{S_\rho\}_{\rho \in I_i}$ are leaves of a local foliation of
some annular neighborhood of $K$.

\medskip

That the construction fails for certain values of $\rho$ is
related to a bifurcation phenomenon. When $k=1$ the families of
surfaces which bifurcate off are (perturbations of) Delaunay
unduloids \cite{Mah}; however, when $k \geq 2$, this bifurcation
is only known to exist in special cases, and the geometry of
the surfaces in the putative bifurcating branches is less
clear. In any case, such bifurcations are inherent
to the problem and occur also in \cite{Mal-Mon-1} and in many
other situations. Furthermore, the index of the hypersurfaces
$S_\rho$, $\rho \in I_i$, tends to $+\infty$ as $i \to \infty$. On
the other hand, we prove that the set $I = \cup_i I_i$ is quite
dense near $0$ in the sense that for any $q \geq 2$ there exists a
$c_q >0$ such that
\[
| {\mathcal H}^1 ((0, \rho) \cap I) - \rho |\leq c_q \, \rho^q,
\]
where ${\mathcal H}^{1} $ denotes the $1$-dimensional Hausdorff
measure.

\medskip

One way to describe the behavior of $S_\rho$ as $\rho$ tends to $0$
is to consider the associated area and curvature densities of $S_\rho$
as $\rho$ tends to $0$; these quantities, properly rescaled, are
extremely close to the corresponding quantities for $\bar S_\rho$,
which in turn satisfy
\begin{equation}
\rho^{k-m} \, {\mathcal H}^{m} \mbox{$\llcorner$} \, \bar S_{\rho}
\, \rightharpoonup \, \omega_{m-k} \, {\mathcal H}^{k} \mbox{
$\llcorner$} \, K \label{eq:00}
\end{equation}
and, for all $q \geq 1$,
\begin{equation}
\rho^{k-m+q} \, |A_{\bar S_{\rho}}|^q \, {\mathcal H}^{m}
\mbox{$\llcorner$}  \, \bar S_{\rho} \, \rightharpoonup \,
(m-k)^{q/2} \, \omega_{m-k} \, {\mathcal H}^{k} \mbox{$\llcorner$}
\, K \label{eq:001}
\end{equation}
as $\rho \searrow 0$. Here $|A_S| ^2 : =  \mbox{Tr} ((A_S)^t \,
A_S)$ is the norm squared of the shape operator. From the explicit
estimates in the construction of $S_\rho$ one can deduce that
(\ref{eq:00}) and (\ref{eq:001}) also hold when $\bar S_\rho$ is
replaced by $S_\rho$.

\medskip

One can ask whether (\ref{eq:00}) and
(\ref{eq:001}) hold for any family of constant mean curvature
hypersurfaces which condense along $K$. It turns out that this is
not the case: families of CMC hypersurfaces condensing along a
nondegenerate geodesic which do not satisfy (\ref{eq:00}) are
constructed in \cite{Mah}. 
In another direction, it is plausible that one should be able
to construct families of CMC hypersurfaces which condense
along lower dimensional sets which are still minimal in
an appropriate sense, but with singularities, for example
a Steiner tree with geodesic edges. A simple example
of this is when $S_\rho$ is obtained by homothetically 
rescaling a fixed Delaunay trinoid in $\R^3$. The limit 
then is a union of three rays meeting at a common vertex, 
each ray having an associated density coming from the 
limiting Delaunay necksize on that end; each ray is minimal, 
of course, and the entire configuration is `balanced' in the 
sense that the weighted sum of the vectors along the 
rays vanishes. 

\medskip

Keeping these various phenomena in mind, it is not clear
whether our main result has a suitable converse, or 
whether it is possible to characterize the possible 
condensation sets of such families of CMC hypersurfaces. 
As a weak and tentative step in this direction we make the

\smallskip

\noindent{\bf Conjecture:} {\sl Let $S_j$ be a family of constant
mean curvature hypersurfaces with mean curvature $H_j \nearrow \infty$;
then for $j$ sufficiently large, $S_j$ is homologically trivial.}

\smallskip

The intuition here is simply that if the $S_j$ were indeed condensing
on a lower dimensional (possibly singular) manifold $K$, then $S_j$
should bound a `tubular neighbourhood' $A_j$ of $K$. 
In any case, this circle of ideas merits further study.

\medskip

In the next section we calculate the asymptotic expansion of the
metric on $M$ in Fermi coordinates around $K$; this is applied in
the (quite technical) \S 3 to derive the expansions of various
geometric quantities for the tubes $\bar S_\rho$ and their
perturbations. This is used in \S 4 to obtain the expression for
the mean curvature of the perturbed tubes, which gives us the
equation which must be solved. An iteration scheme is introduced
in \S 5 which allows us to find a preliminary perturbation for
which the error term is much better, and estimates for the gaps in
the spectrum of the linearization are obtained in \S 6; finally,
the existence of the constant mean curvature hypersurfaces
$S_\rho$ is obtained in \S 7.

\section{Expansion of the metric in Fermi coordinates near $K$}

\subsection{Fermi coordinates}

We now introduce Fermi coordinates in a neighborhood of $K$. For a given $p
\in K$, there is a natural splitting
\[
T_pM = T_p K \oplus N_p K .
\]
Choose orthonormal bases $E_a$, $a=n+1, \ldots, m+1$, for $T_p K$,
and $E_i$, $i=1, \ldots, n$, of $N_pK$.

\medskip

\noindent {\bf Notation : } We shall always use the convention
that indices $a,b,c,d,\ldots \in \{n+1, \ldots, m+1\}$,
indices $i,j,k,\ell,\ldots \in \{1, \ldots, n\}$ and indices
$\alpha, \beta, \gamma, \ldots \in \{1, \ldots, m+1\}$.

\medskip

Consider, in a neighborhood of $p$ in $K$, normal geodesic
coordinates
\[
f(y) : = \exp^K_p (y^a\, E_a), \qquad y := (y^{n +1}, \ldots, y^{m+1}),
\]
where $\exp^K$ is the exponential map on $K$ and summation over repeated
indices is understood. This yields the coordinate vector fields
$X_a : = f_* (\del_{y^a})$. For any $E \in T_p K$, the curve
\[
s \longrightarrow \gamma_E(s) := \exp^K_p (sE),
\]
is a geodesic in $K$, so that
\[
\left. \nabla_{X_a} X_b \right|_p \in N_p K.
\]
We define the numbers $\Gamma_{ab}^i$ by
\[
\left. \nabla_{X_a} X_b \right|_p = \Gamma_{ab}^i \, E_i.
\]

Now extend the $E_i$ along each $\gamma_E(s)$ so that they are parallel
with respect to the induced connection on the normal bundle $NK$.
This yields an orthonormal frame field $X_i$ for $NK$ in a neighborhood of
$p$ in $K$ which satisfies
\[
\left. \nabla_{X_a} X_i \right|_p \in T_p K,
\]
and hence defines coefficients $\Gamma_{ai}^b$ by
\[
\left. \nabla_{X_a} X_i \right|_p = \Gamma_{ai}^b \, E_b.
\]

A coordinate system in a neighborhood of $p$ in $M$ is now defined by
\[
F(x,y) := \exp^M_{f(y)}( x^i \, X_i), \qquad
(x,y) :=(x^1, \ldots, x^n, y^{n+1}, \ldots, y^{m+1}),
\]
with corresponding coordinate vector fields
\[
X_i : = F_* (\del_{x^i}) \qquad \mbox{and} \qquad  X_a : = F_*
(\del_{y^a}).
\]
By construction, $X_\alpha \ |_p =E_\alpha$.

\subsection{Taylor expansion of the metric}

As usual, the Fermi coordinates above are defined so that the
metric coefficients
\[
g_{\alpha \beta} = g( X_\alpha , X_\beta)
\]
equal $\delta_{\alpha\beta}$ at $p$; furthermore,
$g(X_a, X_i) =0$ in some neighborhood of $p$ in $K$.
This implies that
\[
X_b \, g(X_a, X_i) = g(\nabla_{X_b} X_a, X_i) + g(X_a,
\nabla_{X_b}\, X_i) =0
\]
on $K$, which yields the identity
\begin{equation}
\Gamma_{ai}^b = - \Gamma_{ab}^i \label{eq:2.444}
\end{equation}
at $p$.

\medskip

Denote by $\Gamma_a^b : N_p K \longrightarrow {\R}$ the linear form
\[
\Gamma_a^b (E_i) : = \Gamma_{ai}^b
\]
We now compute higher terms in the Taylor expansions of the
functions $g_{\alpha \beta}$. The metric coefficients
at $q := F(x,0)$ are given in terms of geometric data at
$p : = F(0,0)$ and $|x| = \mbox{dist}_g(p,q)$.

\medskip

\noindent {\bf Notation} The symbol $\calO(|x|^r)$ indicates
a function such that it and its partial derivatives of any order,
with respect to the vector fields $X_a$ and $x^i \, X_j$, are
bounded by $ c \, |x|^r$ in some fixed neighborhood of $0$.

\medskip

We begin with the expansion of the covariant derivative~:
\begin{lemma}
At the point of $q=  F(x,0)$, the following expansions hold
\begin{equation}
\begin{array}{rllll}
\nabla_{X_i } \, X_j & = & {\cal O}(|x|) X_\gamma,
\\[3mm]
\nabla_{X_a } \, X_b  & = &  - \Gamma^b_{a} (E_i) \, X_i + {\cal
O}(|x|)
X_\gamma, \\[3mm]
\nabla_{X_a } \, X_i  & = &  \nabla_{X_i} X_a = \Gamma^b_{a}(E_i)
\, X_b + {\cal O}(|x|) X_\gamma, \end{array} \label{eq:3-2a}
\end{equation}
 \label{le:2.1}
\end{lemma}
\noindent{\bf Proof:} We have by construction
\[
\nabla_{X_a} X_b = \Gamma_{ab}^i \, X_i + {\cal O} (|x|) \,
X_\gamma
\]
and
\[
\nabla_{X_a} X_j = \nabla_{X_j} X_a  = \Gamma_{aj}^b \, X_b +
{\cal O} (|x|) \, X_\gamma .
\]
Observe that, because we are using coordinate vector fields,
$\nabla_{X_\alpha}X_\beta = \nabla_{X_\beta}X_\alpha$ for any
$\alpha,\beta$. We also have $ \nabla_{X_i} X_j \ |_p = 0$ since
any $X \in N_{p} K$ is tangent to the geodesic $s \longrightarrow
\exp_{p}^M (s X)$, and hence
\[
\left. \nabla_{X_i+X_j}(X_i + X_j) \right|_p=0
\]
Therefore
\[
\left. (\nabla_{X_i}X_j + \nabla_{X_j}X_i ) \right|_p = 0
\]
This completes the proof of the result. \hfill $\Box$

\medskip

We now give the expansion of the metric coefficients. The
expansion of the $g_{ij}$, $i,j=1, \ldots, n$, agrees with the
well known expansion for the metric in normal coordinates
\cite{Sch-Yau}, \cite{Lee-Par}, \cite{Will}, but we briefly recall
the proof here for completeness.
\begin{proposition}
At the point $q = F(x,0)$, the following expansions hold
\begin{equation}
\begin{array}{rllll}
g_{ij} & = & \delta_{ij} + \frac{1}{3} \, g ( R (E_k ,E_i)\,
E_\ell , E_j ) \, x^k \, x^\ell + {\cal O} (|x|^3)
\\[3mm]
g_{ai} & = & {\cal O} (|x|^2) \\[3mm]
g_{ab} & = & \delta_{ab} + 2 \, \Gamma_{a}^b (E_i) \, x^i + (g( R
(E_k , E_a ) \, E_\ell , E_b ) + \Gamma^c_{a} (E_k) \, \Gamma^b_c
(E_\ell) ) \, x^k \, x^\ell + {\cal O} (|x|^3).
\end{array}
\label{eq:3-1}
\end{equation}
\label{pr:2.1}
\end{proposition}
\noindent{\bf Proof:} By construction, $g_{\alpha \beta}=\delta_{\alpha\beta}$ at $p$, and so
\[
g_{\alpha \beta} =\delta_{\alpha\beta} + {\cal O} (|x|).
\]
Now, from
\[
X_i \, g_{\alpha \beta} = g(\nabla_{X_i} X_\alpha , X_\beta) +
g(X_\alpha, \nabla_{X_i} X_\beta),
\]
Lemma~\ref{le:2.1} and (\ref{eq:2.444}), we get
\[
\left. X_i \, g_{aj} \right|_p =0, \qquad  \left. X_i \, g_{jk} \right|_p =0 \qquad
\mbox{and} \qquad \left. X_i \, g_{a b} \right|_p = \Gamma_{ai}^b +
\Gamma_{ib}^a = 2 \Gamma_{ai}^b.
\]
This yields the first order Taylor expansion
\[
g_{aj} = {\cal O} (|x|^2) , \qquad  g_{ij} = \delta_{ij} + {\cal
O}(|x|^2) \qquad \mbox{and} \qquad g_{ab} = \delta_{ab} + 2 \,
\Gamma_{ai}^b \, x^i + {\cal O} (|x|^2).
\]

To compute the second order terms, it suffices to compute $X_k \,
X_k \, g_{\alpha\beta}$ at $p$ and polarize (i.e.\ replace $X_k$
by $X_i + X_j$, etc.). We compute
\begin{equation}
X_k \, X_k \, g_{\alpha\beta} = g (\nabla_{X_k}^2
X_\alpha,X_\beta) + g(  X_\alpha, \nabla_{X_k}^2 X_\beta) + 2 \,
g( \nabla_{X_k} X_\alpha, \nabla_{X_k} X_\beta) \label{eq:ll}
\end{equation}

To proceed, first observe that
\[
\left. \nabla_{X} X \right|_{p'}= \left. \nabla_X^2 X \right|_{p'}=0
\]
at $p' \in K$, for any $X \in N_{p'} K$. Indeed, for all $p' \in
K$, $X \in N_{p'} K$ is tangent to the geodesic $s \longrightarrow
\exp_{p'}^M(s X)$, and so $\nabla_{X} X = \nabla_X^2 X = 0$ at the
point $p'$.

\medskip

In particular, taking $X = X_k + \e \, X_j$, we obtain
\[
0 = \nabla_{X_k + \e X_j}\nabla_{X_k + \e X_j}(X_k + \e X_j) \
_{|p}
\]
equating the coefficient of $\e$ to $0$ gives
$\nabla_{X_j}\nabla_{X_k}X_k \ _{|p} = -2
\nabla_{X_k}\nabla_{X_k}X_j \ _{|p}$, and hence
\[
\left. 3 \, \nabla_{X_k}^2 X_j \right|_p =  R(E_k,E_j) \, E_k,
\]
So finally, using (\ref{eq:ll}) together with the result of
Lemma~\ref{le:2.1}, we get
\[
\left. X_k \, X_k \, g_{ij} \right|_p = \frac{2}{3} \, g( R(E_k, E_i) \,
E_k , E_j).
\]
The formula for the second order Taylor coefficient for $g_{ij}$
now follows at once.

\medskip

Recall that, since $X_\gamma$ are coordinate vector fields, we
have from (\ref{eq:ll})
\[
\nabla_{X_k}^2 X_\gamma = \nabla_{X_k}\nabla_{X_\gamma }X_k =
\nabla_{X_\gamma} \nabla_{X_k}X_k + R(X_k,X_\gamma) \, X_k.
\]
Using (\ref{eq:ll}), this yields
\[
\begin{array}{rllll}
X_k \, X_k \, g_{ab}  & = & 2  \, g(R(X_k, X_a) X_k, X_b) + 2 \,
g( \nabla_{X_k} X_a , \nabla_{X_k} X_b )
\\[3mm]
&  +  & g (\nabla_{X_a} \nabla_{X_k} X_k ,X_b) + g( X_a ,
\nabla_{X_b } \nabla_{X_k} X_k)
\end{array}
\]
Using the result of Lemma~\ref{le:2.1} together with the fact that
$\nabla_X X =0 \ _{|p'}$ at $p' \in K$  for any $X \in N_{p'} K$,
we conclude that
\[
\left. X_k \, X_k \, g_{ab} \right|_p = 2  \, g(R(E_k, E_a) E_k, E_b) + 2
\, \Gamma_{ak}^c \, \Gamma_{bk}^c
\]
and this gives the formula for the second order Taylor expansion
for $g_{ab}$. \hfill $\Box$

\medskip

Later on, we will need an expansion of some covariant derivatives
which is more accurate than the one given in Lemma~\ref{le:2.1}.
These are given in the~:
\begin{lemma}
At the point $q = F(x,0)$, the following expansion holds
\begin{equation}
\begin{array}{rllll}
 \nabla_{X_a} \, X_b & = & - \Gamma_{a}^b(E_j) \, X_j -  g( R( E_i, E_a) \, E_j , E_b )\, x^i \, X_j \\[3mm]
 & + & \frac{1}{2} \, \left( g ( R (E_a, E_b)\, E_i, E_j) - \Gamma^c_{a} (E_i) \,\Gamma_{c}^b (E_j) - \Gamma^c_{a}(E_j) \,
\Gamma_{c}^b (E_i) \right) \, x^i \, X_j   \\[3mm]
& + & {\cal O} (|x|) \, X_c + {\cal O} (|x|^2) \, X_j.
\end{array}
\label{eq:3-2b}
\end{equation}
\label{le:2.2}
\end{lemma}
\noindent{\bf Proof:} We compute
\[
\begin{array}{rlllll}
X_i \, g( \nabla_{X_a} X_b, X_j ) & = & g( \nabla_{X_i}
\nabla_{X_a} X_b, X_j ) + g( \nabla_{X_a} X_b, \nabla_{X_i} X_j ) \\[3mm]
& = & g (R(X_i, X_a)\, X_b, X_j ) + g( \nabla_{X_a} \nabla_{X_b}
X_i, X_j ) + g( \nabla_{X_a} X_b, \nabla_{X_i} X_j )
\end{array}
\]
Observe that, by construction, we have arranged in such a way that
\[
\nabla_{X_a+ \e X_b} X_i = (\Gamma^c_{ai} + \e \, \Gamma^c_{bi}
)\, X_c
\]
along the geodesic $s \longrightarrow \exp^K_p (s (E_a+ \e E_b))$.
Hence
\begin{equation}
\nabla_{X_a+ \e X_b}^2 X_i = \left( (X_a+ \e \, X_b)(\Gamma^c_{ai}
+ \e \, \Gamma^c_{bi} ) \right) \, X_c + (\Gamma^c_{ai} + \e \,
\Gamma^c_{bi} )\, \nabla_{X_a+ \e \, X_b} \, X_c \label{eq:??}
\end{equation}
Evaluating this at the point $p$ and looking for the coefficient
of $\e$ , we obtain
\[
\left. ( \nabla_{X_a} \, \nabla_{X_b} X_i + \nabla_{X_b} \, \nabla_{X_a}
X_i ) \right|_p - \left.( \Gamma^c_{ai}  \, \nabla_{X_b} \, X_c +
\Gamma^c_{bi}  \, \nabla_{X_a} \, X_c ) \right|_p \in T_p K
\]
Hence we get
\[
\begin{array}{rllll}
\left. g(\nabla_{X_a} \, \nabla_{X_b} X_i, X_j) \right|_p + \left. g(\nabla_{X_b}
\, \nabla_{X_a} X_i, X_j) \right|_p & = & \Gamma^c_{ai} \,
\left. g(\nabla_{X_b} \, X_c, X_j) \right|_p \\[3mm]
& + &  \left. \Gamma^c_{bi}  \, g
(\nabla_{X_a} \, X_c, X_j) \right|_p \\[3mm]
& = & \Gamma^c_{ai} \,\Gamma^j_{bc} + \Gamma^c_{bi}  \,
\Gamma^j_{ac}
\end{array}
\]
Finally, we use the fact that
\[
g(\nabla_{X_b} \, \nabla_{X_a} X_i, X_j) = g ( R(X_b,X_a) \, X_i,
X_j ) +  g ( \nabla_{X_a} \, \nabla_{X_b} X_i , X_j )
\]
to conclude that, at the point $p$
\[
\left. 2 \, g(\nabla_{E_a} \, \nabla_{E_b} E_i, E_j) \right|_p  = g ( R
(E_a, E_b)\, E_i, E_j) +  \Gamma^c_{ai}  \,\Gamma_{bc}^j +
\Gamma^c_{bi} \, \Gamma_{ac}^j
\]
Collecting these estimates together with the fact that
$\left. \nabla_{E_i}E_j \right|_p=0$ we conclude that
\[
\left. 2 \, X_i \, g( \nabla_{X_a} X_b, X_j ) \right|_p = - 2 g (R(E_i,
E_a)\, E_j, E_b ) + g ( R (E_a, E_b)\, E_i, E_j) +  \Gamma^c_{ai}
\,\Gamma_{bc}^j + \Gamma^c_{bi} \, \Gamma_{ac}^j
\]
This easily implies (\ref{eq:3-2b}). \hfill $\Box$

\section{Geometry of tubes}

We derive expansions as $\rho$ tends to $0$ for the metric, second
fundamental form and mean curvature of the tubes $\bar S_\rho$ and
their perturbations.  This is an extension of the computation
in \cite{Maz-Pac}.

\subsection{Perturbed tubes}

We now describe a suitable class of deformations of the geodesic
tubes $\bar S_\rho $, depending on a section $\Phi$ of $NK$ and a
scalar function $w$ on the spherical normal bundle $SNK$.

\medskip

Fix $\rho > 0$. It will be convenient to introduce the scaled variable
$\bar{y} = y/\rho$; we also use a local parametrization $z \rightarrow
\Theta (z)$ of $S^{n-1}$. Now define the map
\[
G(z,\bar y) :=  F \, \big(\rho \,(1+ w(z,\bar y)) \,
\Theta (z) + \Phi(\rho \, \bar y) ,\rho \,\bar y \big),
\]
and denote its image by $S_\rho(w,\Phi)$, so in particular
\[
S_{\rho}(0,0) = \bar S_\rho.
\]

\noindent {\bf Notation : } Because of the definition of these
hypersurfaces using the exponential map, various vector fields
we shall use may be regarded either as fields along $K$ or
along $S_\rho(w,\Phi)$. To help allay this confusion, we
write
\[
\Phi  : = \Phi^j \, E_j \qquad \qquad  \Phi_a : = \partial_{y^a}
\, \Phi^j \, E_j \qquad \qquad \Phi_{ab} : = \partial_{y^a}
\partial_{y^b} \, \Phi^j \, E_j
\]
\[
\Theta  : = \Theta^j \, E_j  \qquad \qquad \Theta_i  :=
\partial_{z^i} \Theta^j \, E_j.
\]
These are all vectors in the tangent space $T_p M$ at the fixed
point $p \in K$.  On the other hand, the vectors
\[
\Psi  : = \Phi^j \, X_j \qquad \qquad  \Psi_a : = \partial_{y^a}
\, \Phi^j \, X_j,
\]
\[
\Upsilon : = \Theta^j \, X_j  \qquad \qquad \Upsilon_i : =
\partial_{z^i} \Theta^j \, X_j
\]
lie in the tangent space $T_q M$, $q = F(z,y)$.

For brevity, we also write
\[
w_j :=  \del_{z^j} w , \quad  w_{\bar a} :=  \del_{\bar y^a} w ,
\qquad w_{ij} : =  \del_{z^i}\, \del_{z^j} w, \quad  w_{\bar a
\bar b} : =  \del_{\bar y^a}\, \del_{\bar y^b} w , \quad w_{\bar a
j} : = \del_{\bar y^a}\, \del_{z^j} w.
\]

In terms of all this notation, the tangent space to $S_\rho(w, \Phi)$
at any point is spanned by the vectors
\begin{equation}
\begin{array}{rcccl}
Z_{\bar a} & = & G_* (\del_{\bar y^a}) & = &  \rho \, ( X_a +
w_{\bar a} \, \Upsilon + \Psi_a ), \qquad a = n+1, \ldots, m+1
\\[3mm]
Z_j &= &  G_* (\del_{z^j}) & = & \rho  \, ((1+w)\, \Upsilon_j +
w_j \, \Upsilon ), \qquad j=1, \ldots, n-1.
\end{array}
\label{eq:defz0zj}
\end{equation}

\subsection{Notation for error terms}
The formulas for the various geometric quantities of $S_\rho(\Phi,w)$
are potentially very complicated, and so it is important to condense
notation as much as possible. Fortunately, we do not need to know
the full structure of all of these quantities. Because it is so
fundamental, we have isolated the notational conventions we shall use
in this separate subsection.

\medskip

Any expression of the form $L(w,\Phi)$ denotes a linear
combination of the functions $w$ together with its derivatives
with respect to the vector fields $\rho \, X_a$ and $X_i$ up to
order $2$, and $\Phi^j$ together with their derivatives with
respect to the vector fields $X_a$ up to order $2$. The
coefficients are assumed to be smooth functions on $SNK$ which are
bounded by a constant independent of $\rho$ in the ${\cal
C}^\infty$ topology (i.e.\ derivatives taken with respect to $X_a$
and $X_i$).

\medskip

Similarly, an expression of the form $Q (w,\Phi)$ denotes a
nonlinear operator in the functions $w$ together with its
derivatives with respect to the vector fields $\rho \, X_a$ and
$X_i$ up to order $2$, and $\Phi^j$ together with their
derivatives with respect to the vector fields $X_a$ up to order
$2$. Again, the coefficients of the Taylor expansion of the
corresponding differential operator are smooth on $SNK$, and $Q$
which vanishes quadratically at $(w,\Phi) = (0,0)$.

\medskip

Finally, any term denoted $\calO(\rho^d)$ is a smooth function on
$SNK$ which is bounded in ${\cal C}^\infty(SNK)$ by a constant times
$\rho^d$.

\subsection{The first fundamental form}

The next step is the computation of the coefficients of the first
fundamental form of $S_\rho (w, \Phi)$. We set
\[
q  : = G(z,0) = F(\rho (1 + w (z,0)) \, \Theta (z) +  \Phi (\rho
z), 0)
\]
and $p := G (0,0)$. We obtain directly from (\ref{eq:3-1}) that
\begin{equation}
\begin{array}{rcl}
g( X_a, X_b)  & = & \delta_{ab} + 2 \, \rho \, \Gamma^b_{a}
(\Theta) + {\cal O} (\rho^2) + 2 \, \Gamma^b_{a} \, (\Phi) +
\rho \, L (w, \Phi) + Q (w, \Phi) \\[3mm]
g ( X_i, X_j) & = &  \delta_{ij} + \frac{\rho^2}{3} \, g( R(\Theta
, E_i) \, \Theta , E_j ) + \calO (\rho^3) \\[3mm]
& + & \frac{\rho}{3} \left( g( R ( \Theta , E_i) \, \Phi ,
E_j ) + g( R(\Phi , E_i ) \, \Theta , E_j ) \right) + \rho^2\, L (w,\Phi) + Q (w,\Phi) \\[3mm]
g( X_i, X_a ) & = &  {\cal O} (\rho^2) + \rho \, L^0(w, \Phi) +
Q^0(w,\Phi).
\end{array}
\label{lems}
\end{equation}

We now explain a simple argument which will be frequently used
throughout the paper. Using the previous expansions, we compute
\[
\begin{array}{rcl}
g (\Upsilon, \Upsilon_j ) & = &  g( \Theta , \Theta_j ) +
\frac{\rho^2}{3} \,  g( R(\Theta , \Theta) \, \Theta ,\Theta_j ) + \calO (\rho^3) \\[3mm]
& + & \frac{\rho}{3} \left( g( R(\Theta,\Theta) \, \Phi ,\Theta_j
) + g ( R(\Phi,\Theta ) \, \Theta , \Theta_j ) \right) + \rho^2\,
L(w,\Phi) + Q(w,\Phi)
\end{array}
\]
However, when $w=0$ and $\Phi=0$, $g (\Upsilon, \Upsilon_j ) =0$
since $\Upsilon$ is normal and $\Upsilon_j$ is tangent to $S_\rho
(0,0)$ then, so that the sum of the first three terms on the
right, which is independent of $w$ and $\Phi$, must also vanish.
This, together with the fact that $R(\Theta , \Theta) =0$ implies
that
\begin{equation}
g( \Upsilon, \Upsilon_j ) = \frac{\rho}{3} \, g ( R(\Phi,\Theta )
\, \Theta , \Theta_j ) + \rho^2\, L (w,\Phi) + Q (w,\Phi)
\label{eq:jjjjj}
\end{equation}

Using similar arguments,  we have
\[
\begin{array}{rcl}
g (\Upsilon, \Upsilon ) & = &  g( \Theta , \Theta ) +
\frac{\rho^2}{3} \,
 g( R(\Theta ,\Theta ) \, \Theta , \Theta_j ) + \calO (\rho^3) \\[3mm]
& + & \frac{\rho}{3} \left( g( R(\Theta , \Theta ) \, \Phi ,
\Theta ) + g ( R(\Phi,\Theta ) \, \Theta ,\Theta ) \right) +
\rho^2\, L (w,\Phi) + Q (w,\Phi)
\end{array}
\]
This, together with the fact that $g(\Upsilon, \Upsilon ) =1$ when
$w=0$ and $\Phi =0$, yields
\begin{equation}
g ( \Upsilon, \Upsilon ) = 1 +  \rho^2 \, L (w, \Phi) + Q (w,
\Phi) \label{eq:edf}
\end{equation}

Using these expansions is is easy to obtain the expansion of the
first fundamental form of $S_\rho(\Phi,w)$.
\begin{proposition}
We have
\begin{equation}
\begin{array}{rlllll}
\ds \rho^{-2} \, g( Z_{\bar a}, Z_{\bar b}) & = & \delta_{ab} + 2
\, \rho \, \Gamma_{a}^b (\Theta ) + \calO(\rho^2) + 2 \,
\Gamma_{a}^b (\Phi ) + \rho \, L(w,\Phi) + Q(w,\Phi)\\[3mm]
\ds \rho^{-2} \, g( Z_{\bar a} , Z_j ) & = & \calO(\rho^2)
+ L(w,\Phi) + Q (w,\Phi) \\[3mm]
\ds \rho^{-2} \,  g ( Z_i, Z_j ) & = & g (\Theta_i , \Theta_j) +
\frac{\rho^2}{3}\, g(R(\Theta , \Theta_i)\, \Theta , \Theta_j)  +
\calO(\rho^3) + 2  \, g (\Theta_i,\Theta_j) \, w  \\[3mm]
& + & \frac{\rho}{3}  \left( g ( R(\Theta , \Theta_i )\Phi ,
\Theta_j ) + g ( R(\Theta , \Theta_j) \Phi , \Theta_i ) \right) +
\rho^2 L (w,\Phi) + Q (w,\Phi).
\end{array}
\label{eq:3-4}
\end{equation}
\label{pr:3-2}
\end{proposition}

\subsection{The normal vector field}

Our next task is to understand the dependence on $(w,\Phi)$
of the unit normal $N$ to $S_\rho (w,\Phi)$.
\begin{proposition}
This unit normal vector field $S_\rho(w, \Phi)$ has the expansion
\begin{equation}
\begin{array}{rrllll}
N &  : = &  - \, \Upsilon +  \alpha^j \, \Upsilon_j + \beta^a \,
X_a + \left( \rho \, L(w,\Phi) + Q(w, \Phi) \right) \, X_a  \\[3mm]
&  + & \left( \rho^2 \, L(w, \Phi) + Q(w,\Phi)\right) \, X_j
\end{array}
\label{eq:3-5}
\end{equation}
where the coefficients $\alpha^j$ are solutions of the system
\[
\alpha^j \, g (\Theta_j, \Theta_i ) = w_i + \frac{\rho}{3} \, g (R
(\Phi, \Theta) \, \Theta, \Theta_i ), \qquad i = 1, \ldots, n-1, \qquad \qquad .
\]
and the coefficients $\beta^a$ are given by
\[
\beta^a  =   w_{\bar a } +  g( \Phi_a , \Theta )
\]
\label{pr:3-3}
\end{proposition}
{\bf Proof~:} Define the vector field
\[
\tilde N  : =  - \, \Upsilon + A^j \, Z_{j} +  B^a \,Z_{\bar a},
\]
and choose the coefficients $A^j$ and $B^a$ so that that $\tilde
N$ is orthogonal to all of the $Z_{\bar b}$ and $Z_i$. This leads
to a linear system for $A^j$ and $B^a$.

\medskip

We have the following expansions
\begin{equation}
\begin{array}{rllllll}
g ( \Upsilon , Z_{\bar a} ) & = & \rho \, w_{\bar a} +
  \rho \, g (\Phi_a , \Theta ) + \rho^2 \,  L(w, \Phi) + \rho \, Q(w, \Phi)\\[3mm]
g ( \Upsilon, Z_j )  & = & \rho \, w_j + \frac{\rho^2}{3} \, g (
R( \Phi  , \Theta ) \, \Theta , \Theta_j ) + \rho^3 \, L(w, \Phi)
+ \rho \, Q (w, \Phi)
\end{array}
\label{eq:dede} \end{equation}
 These follow from (\ref{lems})
together with the fact that $ g(\Upsilon , Z_{\bar a} ) =0$ and $
g( \Upsilon, Z_j ) = 0$ when $w = 0$ and $\Phi = 0$.

\medskip

Using Proposition~\ref{pr:3-2}, we get with little work
\[
B^a = w_{\bar a}  + g( \Theta, \Phi_a) + \rho \, L(w, \Phi) +
\frac{1}{ \rho} \, Q(w, \Phi).
\]
and
\[
A^j \, g( \Theta_j, \Theta_i ) = \frac{1}{ \rho} \, w_i +
\frac{1}{3} \, g( R (\Phi, \Theta) \, \Theta , \Theta_i ) + \rho
\, L(w, \Phi) + \frac{1}{ \rho} \, Q(w, \Phi).
\]
Recall also that $Z_j = \rho \, \Upsilon_j + \rho  \, L(w,\Phi)$
and also that $Z_{\bar a} = \rho \, X_a + \rho \, L(w,\Phi)$.
Collecting these, together with the fact that
\[
|\tilde N|_q  =  1 + \rho^2 \, L(w, \Phi) + Q(w, \Phi).
\]
we obtain
\begin{equation}
\begin{array}{rllll}
N  & : = &  - \, \Upsilon +  \frac{1}{\rho} \, ( \alpha^j \, Z_j +\beta^a \, Z_{\bar a} )
+ \left( L(w,\Phi) +  \frac{1}{\rho} Q (w, \Phi) \right) \, Z_{\bar a}  \\[3mm]
& + &  \left( \rho \, L (w, \Phi) +  \frac{1}{\rho}
Q(w,\Phi)\right) \, Z_j + \left( \rho^2 \, L(w, \Phi) +
Q(w,\Phi)\right) \, \Upsilon
\end{array}
\label{eq:3-5-b}
\end{equation}
The result then follows at once. \hfill $\Box$

\subsection{The second fundamental form}

We now compute the second fundamental form. To simplify the
computations below, we henceforth assume that, at the point
$\Theta (z) \in S^{n-1}$,
\begin{equation}
 g ( \Theta_i, \Theta_j ) =\delta_{ij} \qquad \mbox{and}\qquad
\overline \nabla_{\Theta_i} \Theta_{j} =0 \label{eq:3-55}, \quad
i,j = 1, \ldots, n-1
\end{equation}
(where $\overline \nabla$ is the connection on $TS^{n-1}$).
\begin{proposition}
The following expansions hold
\begin{equation}
\begin{array}{rllll}
\rho^{-2} \,  g( N, \nabla_{Z_{\bar a}} Z_{\bar a} ) & = &
\Gamma_{a}^a (\Theta ) + \rho \, g ( R ( \Theta , E_a) \, \Theta ,
E_a ) + \rho  \, \Gamma_{a}^c (\Theta) \, \Gamma_{c}^a (\Theta) +
{\cal O} (\rho^2) \\[3mm]
& - & \frac{1}{\rho} \, w_{\bar a \bar a}  - g( \Phi_{aa} + R
(\Phi , E_a) \, E_a , \Theta)  +  \Gamma_{a}^c (\Theta) \,
\Gamma_{c}^a (\Phi) - w_j \, \Gamma_{a}^a (\Theta_j) \\[3mm]
& + & \rho \, L(w, \Phi) + \frac{1}{\rho} \, Q(w, \Phi), \\[3mm]

\rho^{-2} \,  g( N, \nabla_{Z_j} Z_{j} ) & = & \frac{1}{\rho} +
\frac{2}{3} \, \rho \,  g( R(\Theta , \Theta_j ) \, \Theta,
\Theta_j ) +  {\cal O} ( \rho^2) \\[3mm]
& - & \frac{1}{\rho}  \, w_{jj}  + \frac{1}{\rho} \, w  +
\frac{2}{3} \, g ( R(\Phi, \Theta_j) \, \Theta , \Theta_j ) \\[3mm]
& + & \rho \, L (w, \Phi)  +  \frac{1}{\rho} \, \, Q (w, \Phi)
\\[3mm]

\rho^{-2} \,  g( N, \nabla_{Z_{\bar a}} Z_{\bar b} ) & = &
\Gamma_{a}^b (\Theta )  - \frac{1}{\rho} \, w_{\bar a \bar b}  + {\cal O} (\rho )
+  L (w, \Phi) + \frac{1}{\rho} \, Q(w, \Phi)  \quad a \neq b  \\[3mm]

\rho^{-2} \, g( N, \nabla_{Z_{\bar a} } Z_j ) & = & {\cal O}(\rho)
+ \frac{1}{\rho} \, L (w, \Phi)  + \frac{1}{\rho} \, Q (w, \Phi) \\[3mm]

\rho^{-2} \,  g( N, \nabla_{Z_i} Z_{j} ) & = & {\cal O} (\rho) +
\frac{1}{\rho}  \, L (w, \Phi)  + \frac{1}{\rho} \, Q(w, \Phi) ,
\quad i \neq j.
\end{array}\label{eq:3-6d}
\end{equation}

\label{pr:3-4}
\end{proposition}
{\bf Proof~:} Some preliminary computations are needed. First note
that by Lemma~\ref{le:2.1}, we have
\begin{equation}
\begin{array}{rllll}
\nabla_{X_a} \, X_b  & = &  - \Gamma^b_{a} (E_i) \, X_i + ( {\cal
O}
(\rho) + L(w, \Phi)+ Q(w, \Phi))  \, X_\gamma, \\[3mm]
\nabla_{X_i} \, X_j  & = &   ( {\cal O}
(\rho) + L(w, \Phi)+ Q(w, \Phi))  \, X_\gamma, \\[3mm]
\nabla_{X_a} \, X_i  & = &  \Gamma^b_{a} (E_i) \, X_b + ( {\cal O}
(\rho) + L (w, \Phi)+ Q(w, \Phi))  \, X_\gamma, \label{eq:3-7}
\end{array}
\end{equation} In particular, this, together with the expression of $Z_{\bar a}$ implies that
\begin{equation}
\begin{array}{lllll}
\nabla_{Z_{\bar a}}  X_i  & = &  \rho \, \Gamma_{a}^b (E_i) \, X_b
+ ({\cal O} (\rho^2) + \rho \, L(w, \Phi)+ \rho \, Q(w, \Phi)) \,
X_\gamma, \\[3mm]
\nabla_{Z_{\bar a}}  X_b  & = &  - \rho \, \Gamma_{a}^b(E_i) \,
X_i + ({\cal O} (\rho^2) + \rho \, L(w, \Phi)+ \rho \, Q(w, \Phi))
\, X_\gamma, \label{eq:3-9}
\end{array}
\end{equation}
We will also need the following expansion which follows from the
result of Lemma~\ref{le:2.2}
\begin{equation}
\begin{array}{rllll}
\nabla_{X_a} \, X_b & = & - \Gamma_{a}^b (Eçj) \, X_j -  g( R( \rho \, \Theta + \Phi , E_a) \, E_j , E_b ) \, X_j \\[3mm]
& + & \frac{1}{2} \, \left( g ( R (E_a, E_b)\, \rho \, \Theta +
\Phi , E_j) - \Gamma^c_{a} (\rho \, \Theta + \Phi) \,\Gamma_{c}^b
(E_j) - \Gamma^b_c(\rho \, \Theta + \Phi) \,
\Gamma_{a}^c (E_j) \right)  \, X_j   \\[3mm]
& + & ( {\cal O} (\rho) + L(w, \Phi) + Q(w, \Phi) ) \, X_c + (
{\cal O} (\rho^2 ) +  \rho \, L(w, \Phi) + Q (w, \Phi) ) \, X_j.
\end{array}
\label{eq:3-2t}
\end{equation}

Finally, we will need the expansions
\begin{equation}
\begin{array}{rllllll}
 g( \Upsilon , X_a) & = & \rho \, L(w, \Phi) +  Q(w, \Phi)\\[3mm]
 g( \Upsilon, \Upsilon_j)  & = & \rho \, L(w, \Phi) + Q (w,
\Phi)
\end{array}
\label{eq:3-10}
\end{equation}
whose proof can be obtained as in \S 3.2, starting from the
estimates  (\ref{lems}).

\medskip

\noindent {\bf First estimate :} We estimate $g( N,
\nabla_{Z_{\bar a}} Z_{\bar b})$ when $a=b$ since the
corresponding estimate, when $a\neq b$ is not as important and
follows from the same proof. We must expand
\[
\rho^{-2} \, g( N, \nabla_{Z_{\bar a}} Z_{\bar a})  = \rho^{-1} \,
\left(  g( N, \nabla_{Z_{\bar a}} X_a ) +  g( N, \nabla_{Z_{\bar
a}} (w_{\bar a} \, \Upsilon) ) +  g( N, \nabla_{Z_{\bar a}} \Psi_a
) \right)
\]
The estimate is broken into three steps:

\medskip

\noindent {\bf Step 1} From (\ref{eq:3-5}), we get
\[
\begin{array}{rlllll}
 g( N, \Upsilon ) & = & -  g( \Upsilon, \Upsilon ) + \alpha^j \, g( \Upsilon_j,
\Upsilon) + \beta^b \, g(X_b, \Upsilon) + (\rho \, L^1 (w, \Phi) +
Q^1 (w,\Phi) ) \, g( X_c ,  \Upsilon ) \\[3mm]
& + &  (\rho^2 \, L(w, \Phi) + Q(w, \Phi) ) \, g( X_j,
\Upsilon )  \\[3mm]
& = & -  1 + \rho^2 \, L(w, \Phi) + Q(w, \Phi)
\end{array}
\]
Substituting $N = -\Upsilon + N + \Upsilon$ gives
\[
g (N, \nabla_{Z_{\bar a}} \Upsilon )  = -  \frac{1}{2} \,
\del_{\bar y^a} g ( \Upsilon, \Upsilon ) + g ( N + \Upsilon,
\nabla_{Z_{\bar a}} \Upsilon )
\]
But it follows from (\ref{eq:edf}) that
\[
\del_{\bar y^a} \, g( \Upsilon, \Upsilon ) = \rho^2 \, L(w, \Phi)+
Q(w, \Phi),
\]
and (\ref{eq:3-9}) together with the expression of $N$ implies
that
\[
g( N + \Upsilon , \nabla_{Z_{\bar a}} \Upsilon ) = \rho \, L(w,
\Phi) + \rho \, Q(w, \Phi).
\]
Collecting these estimates we get
\[
g( N, \nabla_{Z_{\bar a}} \Upsilon )  = \rho \, L(w, \Phi) + Q(w,
\Phi).
\]
Hence we conclude that
\[
g( N, \nabla_{Z_{\bar a}} ( w_{\bar a} \, \Upsilon) ) = w_{\bar a
\bar a} \, g( N, \Upsilon ) + w_{\bar a} \, g ( N, \nabla_{Z_{\bar
a}} \Upsilon ) = - w_{\bar a \bar a} + Q (w, \Phi)
\]

\noindent {\bf Step 2} Next,
\[
g ( N, \nabla_{Z_{\bar a}} \Psi_a ) = \rho \,  g( N , \Psi_{aa} )
+  \Phi^j_a \, g( N, \nabla_{Z_{\bar a}} \, X_j )
\]
From (\ref{eq:3-9}), we have
\[
\Phi^j_a \, g (N, \nabla_{Z_{\bar a}} X_j) = \rho \, L(w, \Phi) +
Q(w,\Phi).
\]
Also, using the decomposition of $N$ and (\ref{lems}), we have
\[
\begin{array}{rllll}
g ( N , \Psi_{aa} ) & = & -  g( \Upsilon , \Psi_{aa} ) +
g ( N + \Upsilon, \Psi_{aa}) \\[3mm]
& = & -  g (\Theta , \Phi_{aa}) + \rho \, L (w, \Phi) +  Q(w,
\Phi))
\end{array}
\]
Collecting these gives
\[
g( N, \nabla_{Z_{\bar a}} \Psi_a )  = - \rho \,  g(\Phi_{aa} ,
\Theta ) + \rho^2 \, L  (w, \Phi) + \rho \, Q (w, \Phi)).
\]

\noindent {\bf Step 3} Expanding $Z_{\bar a}$ gives
\begin{equation}
 g( N, \nabla_{Z_{\bar a}} X_a ) = \rho \,  g( N, \nabla_{X_a} X_a ) + \rho \, w_{\bar a} \, g(N,
\nabla_{\Upsilon} X_a ) + \rho \, \Phi^j_a \, g( N, \nabla_{X_j}
X_b )
\label{eq:3-11}
\end{equation}
With the help of (\ref{eq:3-9}), we evaluate
\[
\begin{array}{rllll}
g( N, \nabla_\Upsilon X_a ) & = & {\cal O} (\rho) + L (w, \Phi) + Q (w,\Phi)\\[3mm]
g( N, \nabla_{X_j} X_a ) & = & {\cal O} (\rho) + L (w, \Phi) + Q (w,\Phi)\\[3mm]
g( N + \Upsilon , \nabla_{X_a} X_a ) & = & - \alpha^j \,
\Gamma_{a}^a \, (\Theta_j) + \rho \, L (w, \Phi) + Q(w,\Phi),
\end{array}
\]
and plugging these into (\ref{eq:3-11}) already gives
\[
g( N, \nabla_{Z_{\bar a}} X_a ) = - \rho \, g( \Upsilon ,
\nabla_{X_a} X_a ) + \rho^2 \, L (w, \Phi) + \rho \, Q (w,\Phi)
\]
Using (\ref{eq:3-2t}) we get the expansion
\[
\begin{array}{rllll}
\nabla_{X_a} X_a  & = & \ds \ds - \Gamma_{a}^a (E_j) \, X_j - g(
R( \rho \Theta + \Phi , E_a) \, E_j , E_a) \, X_j + \Gamma_{a}^c
(\rho \, \Theta + \Phi) \, \Gamma_{c}^a (E_j)\,  X_j \\[3mm]
& + & ({\cal O} (\rho) + L (w, \Phi) + Q (w, \Phi) ) \, X_c +
({\cal O} (\rho^2) + \rho \, L (w, \Phi)  + Q(w, \Phi) ) \, X_j,
\end{array}
\]
Finally, using (\ref{lems}) again, we conclude that
\[
\begin{array}{rlllll}
g( N, \nabla_{Z_{\bar a}} X_b) & = & \rho \, \Gamma_{a}^b (\Theta)
+ \rho^2 \, g(R(\Theta , E_a ) \, \Theta  , E_a ) + {\cal O} (\rho^3) \\[3mm]
& + & \rho \, g(R(\Phi , E_a ) \, \Theta , E_a ) + \rho \,
\Gamma_{a}^c (\rho \, \Theta + \Phi) \, \Gamma_{c}^a (\Theta) -
\rho \, \alpha^j  \, \Gamma_{a}^b  \, (\Theta_j)  \\[3mm]
& + & \rho^2 \, L(w, \Phi) + \rho \, Q(w, \Phi),
\end{array}
\]
which, together with the results of Step 1 and Step 2, completes
the proof of the first estimate.

\medskip

\noindent {\bf Second estimate :} We estimate $g( N, \nabla_{Z_i}
Z_j)$ when $i=j$ since, just as before, the corresponding estimate, when $i\neq j$
is not as important and follows similarly. This part is
taken directly from \cite{Maz-Pac}. Observe that, by
Proposition~\ref{pr:3-2}, we can also write
\[
N = - \Upsilon + \frac{1}{\rho} \, \alpha^j \, Z_j + \hat N,
\]
where
\begin{equation}
\hat N = (L(w, \Phi) + Q(w, \Phi)) \, X_a + (\rho^2 \, L(w, \Phi)
+ Q(w, \Phi) ) \, X_j. \label{eq:hatn}
\end{equation}
Now write
\[
\begin{array}{rlllll}
g( N, \nabla_{Z_j} Z_{j} ) & = &  g( N , \nabla_{Z_{j}} Z_{j} ) \\[3mm]
& =&  g ( \nabla_{Z_j}  \Upsilon , Z_j ) - g( \nabla_{Z_j}
(\alpha^i \, Z_i) , Z_j ) \\[3mm]
& + & g ( \hat N, \nabla_{Z_j} Z_{j}) -\del_{z_j} \, g( \hat N,
Z_{j})
\end{array}
\]

\noindent {\bf Step 1 : } By (\ref{eq:3-7}), we can estimate
\[
\begin{array}{rlllll}
\nabla_{Z_j} Z_{j} & = & \rho \, w_j \, Y_{j} + \rho \, w_{jj} \,
\Upsilon  + \rho \, (1+w) \, \nabla_{Z_j} Y_{j} + \rho \,
w_{j} \, \nabla_{Z_j} \Upsilon\\[3mm]
& = &  ({\cal O} (\rho^3) + \rho^2 \, L(w, \Phi) +
\rho^2 \, L (w, \Phi) \, (L(w, \Phi) + Q(w, \Phi)) ) \, X_a  \\[3mm]
& + & ({\cal O} (\rho^3) + \rho \, L(w, \Phi) + \rho^2 \, Q (w,
\Phi))  \, X_k,
\end{array}
\]
Observe that the coefficient of $X_a$ is slightly better than the
coefficient of $X_k$ since the first two terms only involve the
$X_k$. Using this together with (\ref{eq:hatn}) we conclude that
\[
g( \hat N, \nabla_{Z_j} Z_{j} )  =  \rho^3 \, L(w, \Phi) +
 \rho \, Q(w, \Phi) .
\]

\noindent {\bf Step 2 : } Next, using (\ref{eq:hatn}) together with
(\ref{lems}), we find that
\[
\del_{z_j} \, g( \hat N, Z_{j})  = \rho^3 \, L(w, \Phi) + \rho \,
Q(w, \Phi).
\]

\noindent
{\bf Step 3 : } We now estimate
\[
C :=  2 \,  g( \nabla_{Z_j}  \Upsilon , Z_{j} ) .
\]
It is convenient to define
\[
C' : = \frac{2}{1+w}\,  g ( \nabla_{Z_j} (1+w) \, \Upsilon , Z_{j}
)  ,
\]
It follows from (\ref{eq:dede}) that
\[
C = C' + \rho \, Q (w,\Phi)
\]
hence it is enough to focuss on the estimate of $C'$. To analyze
this term, let us revert for the moment and regard $w$ and $\Phi$
as functions of the coordinates $(z,\bar y)$ and also consider
$\rho$ as a variable instead of just a parameter. Thus we consider
\[
\tilde{F}(\rho,z,\bar y) =  F \big(\rho (1+w(z,\bar y))\Upsilon(
z) + \Phi(t \, \bar y) , t \, \bar y \big).
\]
The coordinate vector fields $Z_j$ are still equal to $\tilde F_*
(\del_{z_j})$, but now we also have $ (1+w)\Upsilon =  \tilde F_*
(\del_\rho)$, which is the identity we wish to use below. Now,
from (\ref{eq:3-4}), we write
\[
C' = \frac{1}{1+w} \, g( \nabla_{\del_\rho} Z_j , Z_{j} ) =
\frac{1}{1+w} \, \del_\rho  g ( Z_j, Z_{j} )
\]
Therefore, it follows from (\ref{eq:3-4}) in
Proposition~\ref{pr:3-2} that
\[
\begin{array}{rlllll}
C  & = &  \frac{1}{1+w} \, \del_\rho \, [ \rho^2 \, g ( \Theta_j,
\Theta_j ) + \frac{\rho^4}{3} \, g ( R(\Theta ,
\Theta_j) \, \Theta, \Theta_j ) + {\cal O} (\rho^5) \\[3mm]
& + & 2 \, \rho^2 \,w \,  g( \Theta_j, \Theta_j ) + \ds
\frac{\rho^3}{3}   \, ( g( R(\Theta , \Theta_j)\, \Phi ,
\Theta_j ) + g ( R(\Theta , \Theta_j) \, \Phi , \Theta_j ) ) \\[3mm]
& + & \rho^4 \, L(w, \Phi) + \rho^2 \, Q(w, \Phi) ] + \rho \, Q(w, \Phi) \\[3mm]
& = & \frac{1}{1+w} \,[ 2 \, \rho \,   g( \Theta_j, \Theta_j ) +
\frac{4}{3} \, \rho^3 \, g ( R(\Theta, \Theta_j) \, \Theta ,
\Theta_j ) + {\cal O} (\rho^4) \\[3mm]
& + & \, 4 \, \rho \,w  \, g ( \Theta_j, \Theta_j )  + \rho^2 \,
\left(  g (R(\Theta, \Theta_j)\, \Phi, \Theta_j ) + g ( R(\Theta
, \Theta_j) \, \Phi, \Theta_j ) \right) \\[3mm]
& + & \rho^3 \, L(w, \Phi) + \rho \, Q(w, \Phi) ] \\[3mm]
& = & \ds 2 \, \rho \, g( \Theta_j, \Theta_j )  + \frac{4}{3} \,
\rho^3 \, g( R(\Theta, \Theta_j) \, \Theta , \Theta_j)
+ {\cal O} (\rho^4)\\[3mm]
& + & \ds  2 \, \rho \,w  \, g ( \Theta_j, \Theta_j ) + \rho^2 \,
( g (R(\Theta, \Theta_j) \, \Phi , \Theta_j ) +  g( R(\Theta ,
\Theta_j)\, \Phi , \Theta_j ) ) \\[3mm]
& + & \rho^3 \, L(w, \Phi) + \rho \, Q(w, \Phi)
\end{array}
\]

\noindent {\bf Step 4 : } Finally, we must compute
\[
\begin{array}{rllll}
D & : = &  2 \, g(  \nabla_{Z_j} (\alpha^i \, Z_i) , Z_j )  \\[3mm]
& =  & 2 \, g ( Z_i , Z_j ) \, \del_{z^j} \alpha^i + 2 \, \alpha^i \,
g( \nabla_{Z_i} Z_j , Z_j) \\[3mm]
& =  &  2 \, g ( Z_i , Z_j )  \, \del_{z^j} \alpha^i + \alpha^i \,
\del_{z^i} \, g ( Z_j , Z_j )
\end{array}
\]
Observe that (\ref{eq:3-55}) implies
\[
\del_{z^j} g ( \Theta_i, \Theta_{j'})  = 0
\]
at the point $p$. Using this together with (\ref{eq:3-4}) and the
expression for the $\alpha^i$ given in Proposition~\ref{pr:3-3},
we get
\[
\alpha^i \, \del_{z^i} \,  g ( Z_j , Z_j )  = \rho^4 \, L(w, \Phi)
+ \rho^2 \, Q(w, \Phi)
\]

It follows from (\ref{eq:3-4}) and the definition of $\alpha^i$
again that
\[
g ( Z_i , Z_j ) \, \del_{z^j} \alpha^i = \rho^2 \, g ( \Theta_i ,
\Theta_j ) \, \del_{z^j} \alpha^i  + \rho^4 \, L (w, \Phi) +
\rho^2 \, Q(w, \Phi)
\]
Therefore, it remains to estimate $ g(  \Upsilon_i , \Upsilon_{j'}
) \, \del_{z^j} \alpha^i$. By definition, we have
\[
\alpha^i \, g ( \Theta_i, \Theta_j) =  w_j + \frac{\rho}{3} \, g(
R(\Phi, \Theta) \, \Theta , \Theta_j )
\]
Differentiating with respect to $z^j$ we get
\begin{equation}
\left(  g( \Theta_i, \Theta_j ) \, \del_{z^j} \alpha^i + \alpha^i
\, \del_{z^j}  g( \Theta_i, \Theta_j ) \right) = w_{jj} +
\frac{\rho}{3} \, \del_{z^j} g ( R(\Phi, \Theta) \, \Theta,
\Theta_j ) \label{eq:start}
\end{equation}
Again, it follows from (\ref{eq:3-55}) that $\del_{z^j}  g(
\Theta_i, \Theta_j ) = 0$. Moreover, using (\ref{eq:3-9}), we
first estimate
\[
\nabla_{Z_j} \Upsilon =  \Upsilon_j + {\cal O}(\rho^2) + \rho \,
L(w, \Phi) +  \rho \, Q(w,\Phi);
\]
and, using in addition (\ref{eq:3-55}), we also get
\[
\nabla_{Z_{j}} \Upsilon_j = a \,\Upsilon + {\cal O}(\rho^2) + \rho
\, L(w, \Phi) + \rho \, Q(w,\Phi)
\]
for some $a \in \R$. Reinserting this in (\ref{eq:start}) yields
\[
\begin{array}{rlllll}
 g( \Theta_i, \Theta_j )  \, \del_{z^j} \alpha^i
& = & w_{jj} + \frac{\rho}{3} \, g( R(\Phi, \Theta_j) \, \Theta ,
\Theta_j ) + \frac{\rho}{3} \,  g( R(\Phi, \Theta ) \, \Theta_j ,
\Theta_j ) +\\[3mm]
& + & \rho^3 \, L(w,\Phi) + \rho^2 \, Q (w,\Phi)),
\end{array}
\]
since $R(\Theta, \Theta) =0$.

\medskip

Collecting these estimates, we conclude that
\[
D = \rho^2 \, w_{jj} \, + \frac{\rho^3}{3} \, g ( R(\Phi, \Theta_j
) \, \Theta, \Theta_j )  + \rho^4 \, L (w, \Phi) + \rho^2 \, Q (w,
\Phi)
\]
since $g( R(\Phi, \Theta) \, \Theta_j , \Theta_j) = 0$. With the
estimates of the previous steps, this finishes the proof of the
estimate.

\medskip

\noindent {\bf Third estimate:} Decompose
\[
\frac{1}{\rho} \,  g (N, \nabla_{Z_{\bar a}} Z_j ) =  g( N,
\Upsilon_j) \, w_{\bar a} + g(N, \Upsilon) \,
 w_{\bar a j}  + (1+w) \,  g(N, \nabla_{Z_{\bar a}} \Upsilon_j ) +
w_j \,  g(N, \nabla_{Z_{\bar a }} \Upsilon ).
\]
As above we use the expression of $N$ given in (\ref{eq:3-5}) to
estimate
\[
 g( N, \Upsilon_j) = -  g( \Upsilon, \Upsilon_j )  +
 g( N + \Upsilon , \Upsilon_j ) = L (w, \Phi) + Q (w, \Phi)
\]
Similarly
\[
 g ( N, \Upsilon ) = - 1 + L (w, \Phi) + Q(w, \Phi)
\]
But now, by (\ref{eq:3-9}), we have
\[
 g( N, \nabla_{Z_{\bar a}} \, \Upsilon_j ) = {\cal O}(\rho^2) + \rho
\, L (w, \Phi) + \rho  \, Q (w, \Phi)
\]
and, as already shown in Step 1
\[
g( N, \nabla_{Z_{\bar a}} \, \Upsilon ) = \rho^2 \, L (w, \Phi) +
Q (w, \Phi),
\]
and the proof of the estimate follows directly. \hfill $\Box$

\section{The mean curvature of perturbed tubes}

Collecting the estimates of the last subsection we obtain the
expansion of the mean curvature of the hypersurface $S_\rho(w ,
\Phi)$. In the coordinate system defined in the previous sections,
we get
\begin{equation}
\begin{array}{rlllllll}
\rho \, m \, H (w, \Phi) & = & n-1 + \rho \, \Gamma^a_a (\Theta) +
\left( g( R ( \Theta , \, E_a ) \, \Theta , E_a) +\frac{1}{3} \,
g( R ( \Theta , \, E_i ) \, \Theta , E_i) \right)
\, \rho^2 \\[3mm]
& - & \Gamma^c_a (\Theta) \, \Gamma^a_c (\Theta) \, \rho^2 + {\cal O} (\rho^3) \notag \\[3mm]
& - & \left(  \rho^2 \, \Delta_K w + \Delta_{S^{n-1}} w + (n-1) \,
w \right) + 2 \, \rho \, \Gamma_a^b (\Theta) \, w_{\bar a \bar b} \\[3mm]
& - & \rho \,  g( \Delta_K \Phi +  R(\Phi, E_a) \, E_a ,
\Theta ) - \Gamma_a^c (\Phi) \, \Gamma_c^a (\Theta)  \\[3mm]
& + & \rho^2 \,  L (w, \Phi) + Q(w, \Phi) .
\notag \label{eq:mc}
\end{array}
\end{equation}

We can simplify this rather complicated expression as follows.
First, note that
\[
K \ \mbox{minimal} \Longleftrightarrow \Gamma^a_a = 0.
\]
Next, define
\begin{equation}
{\cal L}_\rho : = - \left(\rho^2 \, \Delta_K + \Delta_{S^{n-1}} + (n-1)\right),
\label{eq:firstmodel}
\end{equation}
as an operator on the spherical normal bundle $SN K$ with
the expression (\ref{eq:firstmodel}) in any local coordinates.
Also, the Jacobi (linearized mean curvature) operator, for $K$
is defined by
\begin{equation}
{\mathfrak J} : = \Delta^N - {\cal R}^N + {\mathcal B}^N,
\label{eq:jacobi}
\end{equation}
cf. \cite{Law}. To explain the terms here, recall that the Levi-Civita connection
for $g$ induces not only the Levi-Civita connection on $K$, but also a connection
$\nabla^N$ on the normal bundle $NK$. The first term here is simply the rough
Laplacian for this connection, i.e.\
\[
\Delta^N = (\nabla^N)^* \nabla^N.
\]
The second term is the contraction (in normal directions)
of the curvature operator for this connection:
\[
{\cal R}^N  : = \left( R(E_i, \cdot )\, E_i \right)^N,
\]
where the $E_i$ are any orthonormal frame for $N_p K$. Finally,
the second fundamental form
\[
B : T_p K \times T_p K \longrightarrow N_p K , \qquad
B(X,Y) : =  \left( \nabla_X Y \right)^N, \quad X,Y\in T_pK,
\]
defines a symmetric operator
\[
{\cal B}^N : =  B^t \circ B;
\]
in terms of the coefficients $ \Gamma_a^b : =  B(E_a, E_b)$,
\[
g({\cal B}^N  \, X , Y ) = \Gamma_a^b (X)\, \Gamma_b^a (Y).
\]
We also use the Ricci tensor
\[
\mbox{Ric}  (X, Y) = -   g (R (X, E_\gamma)  \, Y, E_\gamma), \qquad
X, Y \in T_p M.
\]

In terms of all of this notation, we have the
\begin{proposition}
Let $K$ be a minimal submanifold. Then the mean curvature of ${\cal T}_\rho
(w, \Phi)$ can be expanded as
\begin{equation}
\begin{array}{rlllllll}
\rho \, m \, H (w, \Phi) & = &  (n-1) + \left( \frac{2}{3} \, g(
{\cal R}^N  \Theta ,  \Theta ) -  \frac{1}{3} \,\mbox{Ric} (\Theta, \Theta) -
g ( {\cal B}^N \Theta, \Theta) \right) \, \rho^2 + {\cal O} (\rho^3) \notag \\[3mm]
& + & {\cal L}_\rho \, w + \rho \,  g( {\mathfrak J} \, \Phi,
\Theta) + {\cal O} (\rho^3) \, \nabla^2_K w + \rho^2 \, L (w ,
\Phi) + Q(w, \Phi). \notag
\end{array}
\label{eq:mcb}
\end{equation}
\label{pr:4.1}
\end{proposition}

The equation $\rho \, m \, H= n-1$ can now be written as
\begin{equation}
\begin{array}{rllll}
{\cal L}_\rho  \, w  + \rho \, g({\mathfrak J} \, \Phi, \Theta) &
=   & - \left( \frac{2}{3} \, g( {\mathcal R}^N \, \Theta , \Theta) -
\frac{1}{3} \,\mbox{Ric} (\Theta , \Theta) - g({\cal B}^N \,
\Theta , \Theta)\right) \, \rho^2 + {\cal O} (\rho^3)  \\[3mm]
&  + & {\cal O} (\rho^3) \,\nabla^2_K w +  \rho^2 \, L (w , \Phi)+
 Q (w , \Phi).
\end{array}
\label{ndf}
\end{equation}

\subsection{Decomposition of functions on $SN K$}

Before proceeding, we now state more clearly our notation
for functions on $SN K$.

Let $(\varphi_j, \lambda_j)$ be the eigendata of $\Delta_{S^{n-1}}$,
with eigenfunctions orthonormal and counted with multiplicity.
These individual eigenfunctions do not make sense on all of $SNK$,
but their span is a well-defined subspace $\calS \subset L^2(SNK)$;
thus $v \in \calS$ if its restriction to each fibre of $SNK$
lies in the span of $\{\varphi_1, \ldots, \varphi_n\}$. We denote by
$\Pi$ and $\Pi^\perp$ the $L^2$ orthogonal projections of $L^2(SNK)$
onto $\calS$ and $\calS^\perp$, respectively.

Now, given any function $v \in L^2(SNK)$, we write
\[
\Pi v = g(\Phi,\Theta), \qquad
\Pi^\perp v = \rho w,
\]
so $v = \rho w + g(\Phi,\Theta)$; here $\Phi$ is a section of the
normal bundle $NK$, and the somewhat elaborate notation in
the second summand here reflects the fact that any element
of $\calS$ can be written (locally) as the inner product
of a section of $NK$ and the vector $\Theta$, whose components
are the linear coordinate functions on each $S^{n-1}$.
We shall often identify  this summand with $\Phi$, and thus,
in the following, $w$ and $\Phi$ will always represent the
components of $v$ in $\calS^\perp$ and $\calS$, respectively. Thus
\[
w = \frac{1}{\rho}\Pi^\perp v, \qquad g(\Phi,\Theta) = \Pi v.
\]

Later on we shall further decompose
\begin{equation}
w = w_0 + w_1
\label{eq:0mode}
\end{equation}
where $w_0$ is a function on $K$ and the integral of $w_1$ over
each fibre of $SNK$ vanishes.

Note that ${\mathfrak J}$ preserves $\calS$ and is invertible
since $K$ is a nondegerate minimal submanifold.

\section{Improvement of the approximate solution}

The first important step in solving (\ref{ndf}) is to use an iteration
scheme to find a sequence of approximate solutions $(w^{(i)},\Phi^{(i)})$
for which the estimates for the error term are increasingly small:
\[
\rho  \, m \, H( w^{(i)}, \Phi^{(i)} ) = n-1 + {\cal O}
(\rho^{i+3}).
\]

Letting $(w^{(0)},\Phi^{(0)})=(0,0)$, we define the sequence $(w^{(i+1)}, \Phi^{(i+1)})
\in \calS^\perp \oplus \calS$ inductively as the unique solution to
\begin{equation}
\begin{array}{rllll}
{\cal L}_0  \, w^{(i+1)} + \rho \, g({\mathfrak J} \, \Phi^{(i+1)}
 , \Theta) & =  & - \left( \frac{2}{3} \, g( {\cal R}^N \Theta ,
\Theta) - \frac{1}{3} \,\mbox{Ric} (\Theta ,
\Theta) - g({\cal B}^N \Theta , \Theta)\right) \, \rho^2 + {\cal O} (\rho^3) \\[3mm]
& - & \rho^2 \, \Delta_K w^{(i)} + {\cal O} (\rho^3) \, \nabla^2_K
w^{(i)}  + \rho^2 \, L (w^{(i)}, \Phi^{(i)}) + Q(w^{(i)};
\Phi^{(i)}) .
\end{array}
\label{eq:popo}
\end{equation}
here
\[
{\cal L}_0 : = - \left(\Delta_{S^{n-1}} + (n-1)\right).
\]
This equation becomes simpler when divided into its $\calS^\perp$ and $\calS$
components. Thus using that $\calL_0$ annihilates $\calS$ and
\[
\frac{2}{3} \, g( {\cal R}^N  \Theta , \Theta) - \frac{1}{3}
\,\mbox{Ric} (\Theta , \Theta) - g( {\cal B}^N \Theta , \Theta)) \in \calS^\perp
\]
since it is quadratic in $\Theta$, (\ref{eq:popo}) can be rewritten as the
two separate equations:
\[
\begin{array}{rllll}
{\cal L}_0  \, w^{(i+1)} & =  & - \left( \frac{2}{3} \, g( {\cal
R}^N  \Theta , \Theta) - \frac{1}{3} \,\mbox{Ric} (\Theta ,
\Theta) - g({\cal B}^N \Theta , \Theta)\right) \, \rho^2 + {\cal O} (\rho^3) \\[3mm]
& - &  \rho^2 \, \Delta_K w^{(i)} + {\cal O} (\rho^3) \,
\nabla^2_K w^{(i)} + \rho^2 \, L( w^{(i)},\Phi^{(i)} ) +
Q(w^{(i)}, \Phi^{(i)}),
\end{array}
\]
and
\[
\begin{array}{rllll}
{\mathfrak J} \, \Phi^{(i+1)} & =  &  {\cal O} (\rho^2) + {\cal O
} (\rho^2) \, \nabla^2_K w^{(i)} + \rho \, L( w^{(i)},\Phi^{(i)} )
+ \rho^{-1} \, Q(w^{(i)}, \Phi^{(i)}) .
\end{array}
\]

\medskip

That there is a unique solution now follows directly from the invertibility of
${\mathfrak J}$ on $\calS$ and ${\cal L}_0$ on $\calS^\perp$, so the only
issue is to obtain estimates.

\begin{lemma}
For this sequence $(w^{(i)},\Phi^{(i)})$, we have the estimates
\[
w^{(i)} = {\cal O} (\rho^2) \qquad \qquad \Phi^{(i)} ={\cal O}
(\rho^2),
\]
\[
w^{(i+1)} - w^{(i)} = {\cal O} (\rho^{i+3}) \qquad  \qquad
\Phi^{(i+1)} - \Phi^{(i)} ={\cal O} (\rho^{i+2})
\]
for all $i \geq 1$. \label{le:lklk}
\end{lemma}
{\bf Proof:} The estimates for $(w^{(1)}, \Phi^{(1)})$ are immediate,
and the result for $i > 1$ is proved by a standard induction using the
general structure of the operators $L$ and $Q$. \hfill $\Box$

\medskip

Finally, replacing $(w, \Phi)$ by $(w^{(i)} + w, \Phi^{(i)} + \Phi)$ in (\ref{ndf}),
the equation we must solve becomes
\begin{equation}
\begin{array}{rllll}
\frac{1}{\rho} \, {\cal L}_\rho  \, w  + g({\mathfrak J} \, \Phi,
\Theta) & =  & {\cal O} (\rho^{i+2}) + {\cal O} (\rho^2) \,
\nabla^2_K w + \rho \, \bar L(w,\Phi) + \frac{1}{\rho} \, \bar
Q(w, \Phi).
\end{array}
\label{eq:sqsq}
\end{equation}
This is of course simply the expansion of the equation
\[
m \, H( w^{(i)}+w, \Phi^{(i)}+\Phi) = \frac{n-1}{\rho}.
\]
The linear and nonlinear operators appearing on the right are
different from the ones before, but enjoy similar properties.

\section{Estimating the spectrum of the linearized operators}

We now examine the mapping properties of the linear operator
\begin{equation}
(w, \Phi) \longmapsto \frac{1}{\rho}  \, {\cal L}_\rho \, w +
 g({\mathfrak J} \, \Phi, \Theta) - {\cal O} (\rho^2) \,
\nabla^2_K w - \rho \, \bar L^2(w,\Phi) \label{eq:op}
\end{equation}
which appears in (\ref{eq:sqsq}). This is not precisely the usual Jacobi
operator (applied to the function $\rho \, w + g(\Phi, \Theta)$),
because we are parametrizing this hypersurface as a graph over
$S_\rho (w^{(i)}, \Phi^{(i)})$ using the vector field $- \Upsilon$
rather than the unit normal.

\medskip

To understand the difference between (\ref{eq:op}) and the Jacobi
operator, recall that if $N$ is the unit normal to a hypersurface
$\Sigma$ and $\tilde N$ is any other transverse vector field, then
hypersurfaces which are ${\mathcal C}^2$ close to $\Sigma$ can be
parameterized as either
\[
\Sigma \ni q \mapsto \exp^M_q(w N) \qquad \mbox{or} \qquad \Sigma
\ni q \mapsto \exp^M_q(\tilde{w} \tilde{N}).
\]
The corresponding linearized mean curvature operators
${\mathbb L}_{\Sigma,N}$ and ${\mathbb L}_{\Sigma,\tilde N}$
are related by
\[
{\mathbb L}_{\Sigma, N}  ( g(N, \tilde N) \, w)  + m \, ({\tilde
N^T} H_\Sigma ) \, w = {\mathbb L}_{\Sigma ,\tilde N} w ,
\]
here $\tilde N^T$ is the orthogonal projection of $\tilde N$ onto
$T\Sigma$. Since ${\mathbb L}_{\Sigma, N}$ is self-adjoint with
respect to the usual inner product, we conclude that $L_{\Sigma,\tilde N}$
is self-adjoint with respect to the inner product
\[
\langle v,w \rangle : = \int_{\Sigma} v \, w \, g(N, \tilde N)\, dA_\Sigma.
\]

Now suppose that $\Sigma = S_\rho (w^{(i)}, \Phi^{(i)})$ and
$\tilde{N} = \Upsilon$. From Lemma~\ref{le:lklk} and Proposition~\ref{pr:3-3}
we have
\[
g (N, - \Upsilon) = 1 + {\cal O} (\rho^2).
\]
Furthermore, from Proposition~\ref{pr:3-2} and
Lemma~\ref{le:lklk}, and the fact that $K$ is minimal, the volume
forms of the tubes $S_\rho (w^{(i)}, \Phi^{(i)})$ and $SN K$ are related by
\[
\sqrt{\mbox{det} (g_{S_\rho (w^{(i)}, \Phi^{(i)})}) } = \rho^{k/2}
\, (1+ {\cal O} (\rho^2)) \, \sqrt{\mbox{det} (g_{SN K})};
\]
hence
\begin{equation}
A_\rho : =  g (N, - \Upsilon) \frac{\sqrt{\mbox{det} (g_{S_\rho
(w^{(i)}, \Phi^{(i)})}) }}{\rho^{k/2} \, \sqrt{\mbox{det}
(g_{SN K})}} = 1 + \calO(\rho^2).
\label{eq:NNN}
\end{equation}

Now define
\begin{equation}
\begin{split}
{\mathbb L}_\rho v = {\mathbb L}_{\rho} ( \rho \, w +  g(\Phi , \Theta)) : =
A_\rho  \, \left( \frac{1}{\rho} \, {\cal L}_\rho \,
 w + g({\mathfrak J} \, \Phi, \Theta) + {\cal O} (\rho^2 ) \,
\nabla^2_K w + \rho \, \bar L (w,\Phi) \right) \\
= \left( \frac{1}{\rho} \, {\cal L}_\rho \,
 w + g({\mathfrak J} \, \Phi, \Theta) + {\cal O} (\rho^2 ) \,
\nabla^2_K w + \rho \, \bar L (w,\Phi) \right),
\end{split}
\label{eq:defLrho}
\end{equation}
where the last equality follows from (\ref{eq:NNN}).

\medskip

Finally, multiplying (\ref{eq:sqsq}) by $A_\rho$ gives one further
equivalent form of this equation,
\begin{equation}
{\mathbb L}_\rho  \, v = {\cal O} (\rho^{2+i}) + \frac{1}{\rho} \,
\tilde Q \left( \frac{1}{\rho} \Pi^{\perp} v, \Pi v \right) ,
\label{eq:sqsq-modif}
\end{equation}
where the nonlinear operator on the right has the same properties as before.

\medskip

Associated to ${\mathbb L}_\rho$ is the quadratic form
\[
{\cal Q}_\rho (w, \Phi) : = \int_{SN K} (\rho \, w + g(\Phi,
\Theta) ) \, {\mathbb L}_{\rho} (\rho \, w + g(\Phi, \Theta) ),
\]
and its corresponding polarization, the bilinear form ${\cal C}_\rho$.
We shall study these forms as perturbations of the model forms
\[
{\cal Q}_0 (w, \Phi) : = \int_{SN K} (\rho^2 \, |\nabla_K w|^2
+ |\nabla_{S^{n-1}} w|^2 - (n-1) \, |w|^2) +
\frac{\omega_{n-1}}{n} \, \int_{K}  g( \mathfrak J \Phi, \Phi)
\]
and associated polarization ${\cal C}_0$.

\medskip

To make precise the sense in which ${\cal Q}_0$ and ${\cal Q}_\rho$
are close, define the weighted norm
\[
\|(w,\Phi)\|^2_{H^1_\rho} : =\int_{SN K} (\rho^2 \, |\nabla_K
w|^2 + |\nabla_{S^{n-1}} w|^2  + |w|^2) + \omega_n \, \int_{K} (
|\nabla_K \Phi|^2 + |\Phi|^2)
\]
and also
\[
\|(w,\Phi)\|^2_{L^2} : =\int_{SN K}  |w|^2  + \omega_n \,
\int_{K} |\Phi|^2.
\]

Using (\ref{eq:NNN}) and the properties of $\bar L$, we have
\begin{equation}
\left|{\cal C}_\rho ((w, \Phi), (w', \Phi')) - {\cal C}_0 ((w, \Phi),
(w', \Phi')) \right| \leq c \, \rho \, \|(w,\Phi)\|_{H^1_\rho} \, \|(w',
\Phi')\|_{H^1_\rho} \label{eq:5.55}
\end{equation}

\subsection{Estimates for eigenfunctions with small eigenvalues}
\begin{lemma}
Let $\sigma$ be an eigenvalue of ${\mathbb L}_\rho$ and $(w,\Phi)$
a corresponding eigenfunction. There exist constants $c, c_0 >0$
such that if $|\sigma|\leq c_0$, then using the decomposition $w = w_0 + w_1$
from (\ref{eq:0mode}),
\[
\|(w-w_{0},\Phi)\|_{H^1_\rho}\leq c\,\rho\,\|(w,\Phi)\|_{H^1_\rho}.
\]
\label{le:local}
\end{lemma}
{\bf Proof:} For any $(w',\Phi')$,
\[
\begin{array}{llll}
\calC_\rho((w,\Phi),(w',\Phi')) & = & \displaystyle \sigma
\int_{SN K} (
\rho^2 w \, w'  + g(\Phi,\Theta)g(\Phi',\Theta) ) \\[3mm]
& = & \displaystyle \sigma \int_{SN K} \rho^2 w \, w'  +
 \sigma \, \frac{\omega_n}{n} \, \int_K g(\Phi,\Phi').
\end{array}
\]
In addition, (\ref{eq:5.55}) gives
\begin{equation}
\begin{array}{llll}
\displaystyle \left| \int_{SN K} (\rho^2 \, \nabla_K w \,
\nabla_K w' + \nabla_{S^{n-1}} w  \nabla_{S^{n-1}} w' - (n-1 +
\sigma) \, w \, w') \right. \\[3mm]
\displaystyle \left. \qquad \qquad \qquad \qquad \ds +
\frac{\omega_{n-1}}{n}  \, \int_{K} (g( \mathfrak J \Phi,
\Phi') - \sigma \, g(\Phi, \Phi') ) \right| \leq c \, \rho \,
\|(w, \Phi)\|_{H^1_\rho} \, \|(w',\Phi')\|_{H^1_\rho}.
\end{array}
\label{eq:fdfd}
\end{equation}

\noindent {\bf Step 1 : } Take $w'=0$ and $\Phi' =\Phi$ in
(\ref{eq:fdfd}); this yields
\[
\left| \int_{K} ( g( \mathfrak J \Phi, \Phi) + \sigma \, g(\Phi ,
\Phi) ) \right| \leq c \, \rho  \, \|(w, \Phi)\|_{H^1_\rho} \,
\|(0,\Phi)\|_{H^1_\rho}
\]
Since ${\mathfrak J}$ is invertible, there exists $c_1 >0$ such that
\[
2 \, c_1 \, \|(0,\Phi)\|^2_{H^1_\rho} \leq \left| \int_{K} g(
\mathfrak J \Phi, \Phi) \right|,
\]
hence
\[
(2 \, c_1 - |\sigma |) \, \|(0,\Phi)\|_{H^1_\rho} \leq c \, \rho \,
\|(w,\Phi)\|_{H^1_\rho}
\]
Assuming $c_1 \, \geq  |\sigma|$, we conclude that
\[
\| (0,\Phi) \|_{H^1_\rho} \leq c \, \rho \, \|(w, \Phi)\|_{H^1_\rho}
\]

\noindent {\bf Step 2 : } Now use (\ref{eq:fdfd}) with $\Phi' =0$
and $w = w_1$ to get
\[
\left| \int_{SN K} (\rho^2 \, |\nabla_K w_1|^2 +
|\nabla_{S^{n-1}} w_1|^2 - (n-1-\sigma) \, |w_1|^2  ) \right| \leq
c \, \rho  \, \|(w,\Phi)\|_{H^1_\rho}\,\|(w_1,0)\|_{H^1_\rho}.
\]
However, since $\Pi w_1 = 0$ and $\int_{S^{n-1}}w_1 = 0$, we have
\[
\int_{S^{n-1}} |\nabla_{S^{n-1}} w_1|^2 \geq 2\, n \,
\int_{S^{n-1}} |w_1|^2,
\]
hence
\[
\left| \int_{SN K} (\rho^2 \, |\nabla_K w_1|^2 + \frac{1}{2} \,
|\nabla_{S^{n-1}} w_1|^2 + (1 - |\sigma|) \, |w_1|^2 )\right| \leq
c \,\rho\,\|(w,\Phi)\|_{H^1_\rho} \, \|(w_1,0)\|_{H^1_\rho}.
\]
This implies that
\[
\|(w_1,0)\|_{H^1_\rho}  \leq c \, \rho  \, \|(w,\Phi)\|_{H^1_\rho}
\]
provided $|\sigma | \leq 1/4$. This completes the proof if $c_0
= \mbox{min} (c_1, 1/4)$. \hfill $\Box$

\subsection{Variation of small eigenvalues with respect to $\rho$}

We shall need to obtain some information about the spectral
gaps of ${\mathbb L}_\rho$ when $\rho$ is small, and to do this,
it is necessary to understand the rate of variation of the small eigenvalues
of this operator.
\begin{lemma}
There exist constants $c_0,c >0$ such that, if $\sigma$ is an eigenvalue of
${\mathbb L}_\rho$ with $|\sigma| < c_0$, then
\[
\rho \, \del_\rho \sigma \geq 2 \, (n-1) - c \, \rho
\]
provided $\rho$ is small enough. \label{le:evol}
\end{lemma}
{\bf Proof:} There is a well-known formula for the variation of a simple eigenvalue;
complications arise in the presence of multiplicities, but a result of Kato
\cite{Kato} shows that if one considers the derivative of the eigenvalue as a multi-valued
function, then an analogue of this same formula holds:
\[
\del_\rho  \sigma \in \left\{ \int_{SN K} v \, (\del_\rho
{\mathbb L}_\rho ) v  \quad  : \quad v = \rho \, w + g( \Phi,
\Theta), \qquad {\mathbb L}_\rho v = \sigma \, v, \qquad \| v
\|_{L^2} =1 \right\}.
\]
Hence we must provide bounds for the set on the right. We do this
by comparing to the model case and using the bounds for eigenfunctions
obtained in the last subsection.

\medskip

Let ${\mathbb L}_\rho v = \sigma \, v$, but rather than normalizing
by $\|v \|_{L^2} =1$), assume instead that $\|(w,\Phi)\|_{L^2} =1$.
In order to compute $\del_\rho {\mathbb L}_\rho$, recall that
$w = \rho^{-1}\Pi^\perp v$, so we can write
\[
{\mathbb L}_\rho v = \frac{1}{\rho^2}\calL_\rho \Pi^\perp v +
g({\mathfrak J}\Phi,\Theta) + \calO(\rho)\nabla_K^2 \Pi^\perp v +
\rho \, \bar{L}( \rho^{-1} \, \Pi^\perp v, \Pi v).
\]
Since $\Pi$ and $\Pi^\perp$ are independent of $\rho$, we have
\[
\del_\rho {\mathbb L}_\rho v =
-\frac{2}{\rho^3}\calL_\rho(\Pi^\perp v) + \frac{1}{\rho^2}(-2\rho
\Delta_K \Pi^\perp v) + \calO (1) \, \nabla^2_K \Pi^\perp v + \,
\bar{L}( \rho^{-1} \, \Pi^\perp v,  \Pi v )
\]
\[
= - \frac{2}{\rho^2}\calL_0 w + \calO(\rho)\nabla^2_K w + \bar{L}( w, \Phi).
\]
where the operator $\bar L$ varies from line to line but satisfies
the usual assumptions. This now gives
\begin{equation}
\left| \int_{SN K} v \, (\del_\rho {\mathbb L}_\rho ) v +
\frac{2}{\rho} \, \int_{SN K} (|\nabla_{S^{n-1}} w|^2 - (n-1) \, |w|^2)
\right|\leq c \, \|( w, \Phi ) \|^2_{H^1_\rho}. \label{eq:cc}
\end{equation}

Now, for this eigenfunction $v$, $Q_\rho(v,v) = \sigma \int \rho^2 |w|^2 + g(\Phi,\Phi)$,
and hence by (\ref{eq:5.55}),
\begin{equation}
\begin{split}
\displaystyle \left| \int_{SN K} (\rho^2 \, |\nabla_K w|^2 +
|\nabla_{S^{n-1}} w|^2 - (n-1 + \sigma ) \, |w|^2) +
\frac{\omega_{n-1}}{n} \, \int_{K} (g( \mathfrak J \Phi, \Phi)
-  \sigma \, g(\Phi, \Phi) ) \right| \\
\leq c \, \rho \, \|(w, \Phi)\|_{H^1_\rho}^2,
\end{split}
\label{eq:cccc}
\end{equation}

By Lemma~\ref{le:local},
\begin{equation}
\int_{SN K} |\nabla_{S^{n-1}} w|^2 + \int_K (|\nabla_K \Phi|^2+
|\Phi|^2) \leq c \, \rho \, \|(w,\Phi)\|_{H^1_\rho}^2,
\label{eq:ppp}
\end{equation}
and inserting this in (\ref{eq:cccc}) gives
\begin{equation}
\left| \int_{SN K} (\rho^2 \, |\nabla_K w|^2 - (n-1+\sigma ) \,
|w|^2) \right| \leq  c \, \rho \, \|(w, \Phi)\|_{H^1_\rho}^2.
\label{eq:szsz}
\end{equation}
Adding these last two estimates now implies that
\[
\|(w, \Phi)\|_{H^1_\rho}^2  \leq  c \, \rho \,  \|(w,
\Phi)\|_{H^1_\rho}^2 + c \, \int_{SN K} |w|^2;
\]
Thus, when $\rho$ is small enough,
\[
\|(w, \Phi)\|_{H^1_\rho}^2 \leq c \| (w, \Phi)\|_{L^2} \leq c
\]
by our choice of normalization. From (\ref{eq:ppp}) again
\[
\int_{SN K} |\nabla_{S^{n-1}} w|^2 + \int_K (|\nabla_K \Phi|^2+
|\Phi|^2) \leq c \, \rho.
\]

Inserting this into (\ref{eq:cc}), and using again that $\| (w, \Phi)\|_{L^2} =1$,
we get
\begin{equation}
\left| \int_{SN K} v \, (\del_\rho {\mathbb L}_\rho ) v -
\frac{2}{\rho} \, (n-1) \right|\leq c
\label{frfr}
\end{equation}
for all $v$ such that ${\mathbb L}_\rho v = \sigma \, v$ and $\|(w,\Phi)\|_{L^2} =1$.

\medskip

This already implies that $\del_\rho \sigma >0$ for $\rho$ small
enough. But observing that we always have $||v||_{L^2} \leq  \|(w, \Phi)\|_{L^2}$,
we conclude that
\[
\inf_{ \stackrel{v: {\mathbb L}_\rho v = \sigma}{\|v\|_{L^2} = 1}}
\int_{SN K} v \, (\del_\rho {\mathbb L}_\rho ) v \quad  \geq
\inf_{\stackrel{v: {\mathbb L}_\rho v = \sigma \, v}{\|(w,\Phi)\|_{L^2} = 1}}
\int_{SN K} v \, (\del_\rho {\mathbb L}_\rho) v
\]
and (\ref{frfr}) implies that
\[
\del_\rho \sigma \geq  \frac{2}{\rho} \, (n-1) - c .
\]
This completes the proof of the result. \hfill $\Box$

\subsection{The spectral gap at $0$ of ${\mathbb L}_\rho$}

We can now prove a quantitative statement about the clustering
of the spectrum at $0$ of ${\mathbb L}_\rho$ as $\rho \searrow 0$.
The ultimate goal is to estimate the norm of the inverse of
this operator, but by self-adjointness, this is equivalent to
an estimate on the size of the spectral gap at $0$.

\begin{lemma}
Fix any $q \geq 2$. Then there exists a sequence of disjoint
nonempty intervals $I_i = (\rho_i^-, \rho_i^+)$,  $\rho_i^\pm
\rightarrow 0$ and a constant $c_q > 0$ such that when $\rho \in I
: =  \cup_i I_i$, the operator ${\mathbb L}_\rho$ is invertible
and
\[
{\mathbb L}_\rho^{-1} : L^2(SN K) \longrightarrow L^2 (SN K)
\]
has norm bounded by $c_q \,  \rho^{-k-q+1}$, uniformly in $\rho
\in I$. Furthermore, $I  :=  \cup_i I_i$ satisfies
\[
\left|{\mathcal H}^1 ( (0, \rho) \cap I) -\rho \right|\leq c \,
\rho^{q}, \qquad \rho \searrow 0.
\]
\label{le:sg}
\end{lemma}
{\bf Proof:} An estimate for the size of the spectral gap at $0$
is related to the spectral flow of ${\mathbb L}_\rho$, and so it
suffices to find an asymptotic estimate for the number of negative
eigenvalues of ${\mathbb L}_\rho$. Define the two quadratic forms
\[
{\cal Q}^{\pm} (w, \Phi) : =  {\cal Q}_0 (w, \Phi) \pm  \gamma \, \rho
\, \|(w, \Phi) \|_{H^1_\rho}^2
\]
From (\ref{eq:5.55}), if $\gamma >0$ is sufficiently large, then
\[
{\cal Q}^- \leq {\cal Q}_\rho \leq {\cal Q}^+,
\]
and this will give a two-sided bound for the index of ${\cal Q}_\rho$.

\medskip

Decomposing $w = w_0 + w_1$ with $w_0$ depending only on $y\in K$, we write
\[
D^\pm_0(w_0) = (1 \pm \gamma \, \rho ) \, \int_{K} \rho^2 \,
|\nabla_K w_0|^2  - (n-1 \mp \gamma \, \rho) \, \int_{K} |w_0|^2,
\]
\[
D^\pm_1(w_1) = (1 \pm \gamma \, \rho ) \, \int_{SN K} (\rho^2
\, |\nabla_K w_1|^2 + |\nabla_{S^{n-1}}  w_1|^2 ) - (n-1 \mp
\gamma \, \rho) \, \int_{SN K} |w_1|^2,
\]
and finally
\[
D^\pm(\Phi) = - ( 1\pm \gamma \, \rho) \, \int_K g({\mathfrak J}
\, \Phi, \Phi),
\]
so that
\[
{\cal Q}^\pm (w, \Phi) = \omega_{n-1} \, D^\pm_0 (w_0) +
D^\pm_1(w_1)+ \frac{\omega_{n-1}}{n} \,  D^\pm (\Phi)
\]
If $1 - \gamma \,\rho >0$, then the index of $D^\pm$ equals the
index of the minimal submanifold $K$, and hence does not depend on
$\rho$. Next, if $(1 - \gamma\, \rho) \, 2\, n - (n - 1 + \gamma
\, \rho) >0$, then the index of $D_1^\pm$ equals $0$. So it
remains only to study the index of $D_0^\pm$. This is equal to
the largest $j \in \N$ such that
\[
(1 \pm \gamma \, \rho ) \, \rho^2 \, \mu_j  \leq (n-1 \mp c \,
\rho)
\]
Weyl's asymptotic formula states that
\[
\mbox{Ind} \, {\mathcal Q}^\pm \sim c_K \, \rho^{-k},
\]
and hence the index of $D_0^\pm$, and finally $\mbox{Ind} \, {\cal Q}_\rho$
too, is asymptotic to $c_K \,\rho^{-k}$.

Let $\rho_i \searrow 0$ be the decreasing sequence corresponding
to the values at which the index of ${\cal Q}_\rho$ changes, counted
according to the dimension of the nullspace of ${\mathbb L}_{\rho_i}$, i.e.\
\[
\rho_{i-1} < \rho_{i} = \ldots =  \rho_j < \rho_{j+1}
\]
if $\mbox{dim Ker} \, {\mathbb L}_{\rho_i} = j+1-i$. This is
well-defined since, by Lemma~\ref{le:evol} the small
eigenvalues of ${\mathbb L}_\rho$ are monotone increasing for
$\rho$ small enough and hence, the function $\rho \rightarrow
{\cal Q}_\rho$ is monotone decreasing for $\rho$ small.

\medskip

The estimates for $\mbox{Ind} \, Q_{2\rho}$ and $\mbox{Ind} \,
Q_{\rho}$ imply that
\[
r_\rho := \# \{ \rho_i\in (\rho, 2\rho) \} \sim c\, \rho^{-k}.
\]
Letting $l_\rho$ denote the sum of lengths of intervals
$(\rho_{i+1},\rho_i)$ for which $\rho_{i+1} \in (\rho, 2 \rho)$
and $(\rho_i -\rho_{i+1}) \leq \rho^{k+q}$, then we have $l_\rho
\leq c \, \rho^q$; from this we conclude that $\ell_\sigma$, the
sum of lengths of all intervals $(\rho_{i+1},\rho_i)$ where
$\rho_{i+1} < \rho$ and $(\rho_i - \rho_{i+1}) \leq \rho^{k+q}$ is
also estimated by $c \, \rho^{q}$.

Define
\[
\tilde I = \bigcup_{i \in J} (\rho_{i+1}, \rho_i), \qquad
\mbox{where}\qquad i \in J \Leftrightarrow \rho_i - \rho_{i+1}
\geq \rho_i^{k+q}.
\]
Then by the above, we have
\[
\left|{\mathcal H}^1 ( (0,\rho) \cap I) -\rho \right| \leq c_q \,
\rho^q.
\]

Finally, consider for any $\rho \in (\rho_{i+1}, \rho_i)$, $i \in
J$, the eigenvalues of ${\mathbb L}_\rho$ which are closest to
$0$, say
\[
\sigma^-(\rho) < 0 < \sigma^+(\rho).
\]
(Thus for each $\rho \in (\rho_{i+1},\rho_i)$, $\sigma^-(\rho) =
\sigma_j$ where $j = \mbox{Ind}\,Q_\rho$.) By construction,
\[
\lim_{\rho \searrow \rho_{i+1}} \sigma^+ (\rho) = \lim_{\rho
\nearrow \rho_i} \sigma^- (\rho) = 0.
\]
By Lemma~\ref{le:evol},
\[
\sigma^-(\rho) \leq  2 \, (n-1) \,\frac{\rho-\rho_{i}}{\rho_{i}} +
c \, \rho_{i}^{k+q}, \qquad \rho \in (\rho_{i+1}, \rho_i),
\]
and
\[
\sigma^+(\rho) \geq  2 \, (n-1) \,
\frac{\rho-\rho_{i+1}}{\rho_{i+1}} - c \, \rho_{i+1}^{k+q}, \qquad
\rho \in (\rho_{i+1},\rho_{i}).
\]
Hence by the monotonicity of small eigenvalues, if
\[
\rho \in I : = \bigcup_i (\rho_{i+1} + \frac{1}{4} \,
\rho_i^{k+q}, \rho_i - \frac{1}{4} \, \rho_i^{k+q})
\]
then the infimum of the absolute value of the eigenvalues of
${\mathbb L}_\rho$ is bounded from below by a constant (only
depending on $n$) times $\rho_i^{k+q-1}$, provided $\rho$ is small
enough. The result then follows at once. \hfill $\Box$

\section{Existence of constant mean curvature hypersurfaces}

We now use the results of the previous sections in order to solve
the equation (\ref{eq:sqsq-modif}) which reduces to find a fixed
point
\[
 \rho \, w  + g(\Phi, \Theta)  = {\mathbb L}_\rho^{-1} \left(
{\cal O} (\rho^{2+i}) + \frac{1}{\rho} \, \tilde Q(w, \Phi)
\right).
\]

Since any function $v$ defined on $SN K$ can be decomposed as
$v = \rho \, w + g(\Phi, \Theta)$ where the function $w$ satisfies
\[
\int_{S^{n-1}} w \, \varphi_j =0
\]
for all $j=1, \ldots, n$, this equation can be re-written as
\[
v  =  {\mathbb L}_\rho^{-1} \left( {\cal O} (\rho^{2+i}) +
\frac{1}{\rho} \, \tilde Q \left( \frac{1}{\rho} \, \Pi^\perp v,
\Pi v \right) \right)
\]
We start with the following elementary observation
\begin{lemma}
There exists a constant $c >0$ such that
\[
\rho^{2+\alpha} \, \| v \|_{{\cal C}^{2,\alpha }} \leq c \,
 \rho^2 \, \| {\mathbb L}_\rho \, v \|_{{\cal C}^{0, \alpha}} + c \, \rho^{-\frac{k}{2}} \,
\| v \|_{L^2}
\]
\label{le:vvc}
\end{lemma}
{\bf Proof :} This is a simple application of (rescaled) standard
elliptic estimates. We set $f : = {\mathbb L}_\rho \, v$ and, as
in \S 3.1, we use local normal coordinates $\bar y = y/\rho$ to
parameterize a ball of radius $2 \, \rho \, R$ in $K$, for some
fixed small constant $R >0$, and local coordinates $z$ to
parameterize $S^{n-1}$. Define the functions
\[
\bar v (z, \bar y ) : = v (z, \rho \, \bar y) \qquad \mbox{and}
\qquad \bar f (z, \bar y ) : = \rho^2 \, f (z, \rho \, \bar y)
\]
It is easy to check that $f : = {\mathbb L}_\rho  \, v$ translates
into $\bar {\mathbb L}_\rho \bar v = \bar f$, where $\bar {\mathbb
L}_\rho$ is a second order elliptic operator whose coefficients
are bounded uniformly in $\rho$ as $\rho$ tends to $0$. Moreover,
the principal part of $\bar {\mathbb L}_\rho$ is the Laplace
operator on $SN K$. Standard elliptic estimates yield
\[
\| \bar v \|_{\bar {\cal C}^{2,\alpha }(B_R \times S^{n-1})} \leq
c \, \| \bar f \|_{\bar {\cal C}^{0, \alpha}(B_R \times S^{n-1})}
+ c \, \left(\int_{S^{n-1}} \, \left(\int_{B_{2R}} |\bar v|^2 \,
d\bar y \right) \right)^{1/2}
\]
where, to evaluate the H\"older norms in $\bar {\cal C}^{p,\alpha
}$ one takes derivatives with respect to $\bar y$ and $z$. Going
back to the functions $v$ and $f$ we have
\[
\rho^{2+\alpha} \, \| v \|_{{\cal C}^{2,\alpha }(B_{\rho R} \times
S^{n-1})} \leq c \, \| \bar v \|_{\bar {\cal C}^{2,\alpha }(B_{R}
\times S^{n-1})}, \qquad \| \bar f \|_{\bar {\cal
C}^{2,\alpha }(B_{\rho R} \times S^{n-1})} \leq c \, \rho^2 \, \|
f \|_{{\cal C}^{2,\alpha }(B_{R} \times S^{n-1})}
\]
and
\[
\left(\int_{S^{n-1}} \, \left(\int_{B_{2R}} |\bar v|^2 \, d\bar y
\right) \right)^{1/2} \leq c \, \rho^{-\frac{k}{2}}
\left(\int_{S^{n-1}} \, \left(\int_{B_{2\rho R}} |v|^2 \, dy
\right) \right)^{1/2}
\]
the result then follows at once. \hfill $\Box$

\medskip

We fix $q \geq 2$ and $\alpha \in (0,1)$. Collecting the result of
Lemma~\ref{le:sg} and the result of the previous Lemma, we
conclude that, if $\rho \in I$, then
\begin{equation}
\| v \|_{{\cal C}^{2,\alpha }} \leq c \, \rho^{-D} \, \| {\mathbb
L}_\rho \, v \|_{{\cal C}^{0, \alpha}} \label{estd}
\end{equation}
where the constant $c>0$ does not depend on $\rho$ and where $D :
= 3 \frac{k}{2} + q + 1 + \alpha $.

\medskip

Given $R >0$, set
\[
B (R)  : = \{ v \in  {\cal C}^{2,\alpha} ( SN K)  \, : \, \| v
\|_{{\cal C}^{2, \alpha}}\leq  R \}.
\]
and define the mapping
\[
{\cal N}_\rho (v)  : =  {\mathbb L}_\rho^{-1} \left( {\cal O}
(\rho^{2+i}) + \frac{1}{\rho} \, \tilde Q \left( \frac{1}{\rho} \,
\Pi^\perp v, \Pi v \right) \right)
\]
It follows from (\ref{estd}) that we have
\[
\| {\cal N}_\rho ( 0) \|_{{\cal C}^{2,\alpha }} \leq \frac{c_0}{2}
\, \rho^{2+i-D}
\]
for some constant $c_0 >0$, independent of $\rho \in I$.

\medskip

We choose $i \in {\mathbb N}$ such that $i > 2 \, D +1$. Using the
properties of the operator $\tilde Q$, it is easy to check that
there exists $\rho_0 >0$ such that, for all  $\rho \in (0, \rho_0)
\cap I$,
\[
\| {\cal N}_\rho ( v ) \|_{{\cal C}^{2,\alpha }} \leq c_0 \,
\rho^{2+i-D}
\]
and
\[
\| {\cal N}_\rho ( v ) - {\cal N}_\rho (v') \|_{{\cal C}^{2,\alpha
}} \leq c \, \rho^{i-1-2 D} \, \|v-v'\|_{{\cal C}^{2,\alpha }}
\]
for all $v, v' \in B (c_0\, \rho^{2+i-D})$. Therefore the mapping
${\cal N}_\rho$ admits a (unique) fixed point $v_\rho$ in $B (c_0
\, \rho^{2+i-D})$. This yields the existence of a constant mean
curvature perturbation of the tube $S_\rho (w^{(i)}, \Phi^{(i)})$
for all $\rho \in (0, \rho_0) \cap I$. The proof of the Theorem is
complete.

\end{document}